\documentclass{amsart}

\usepackage{amsmath,amssymb,amsthm}
\usepackage{graphicx,tikz}
\newtheorem{theorem}{Theorem}

\newtheorem*{corollary}{Corollary}

\theoremstyle{definition}

\theoremstyle{remark}

\begin{document}

\title[]{The product of two high-frequency Graph Laplacian eigenfunctions is smooth}
\thanks{S.S. is supported by the NSF (DMS-2123224) and the Alfred P. Sloan Foundation.}

\author[]{Stefan Steinerberger}
\address{Department of Mathematics, University of Washington, Seattle, WA 98195, USA}
\email{steinerb@uw.edu}

\begin{abstract} In the continuous setting, we expect the product of two oscillating functions to oscillate
even more (generically). On a graph $G=(V,E)$, there are only $|V|$ eigenvectors
of the Laplacian $L=D-A$, so one oscillates `the most'. The purpose of this short note is to point out an interesting phenomenon: if $\phi_1, \phi_2$ are delocalized eigenvectors of $L$ corresponding to large eigenvalues, then their (pointwise) product $\phi_1 \cdot \phi_2$ is smooth (in the sense
of small Dirichlet energy): 
highly oscillatory functions have largely matching oscillation patterns.
\end{abstract}

\maketitle

\section{Introduction and Result}
\subsection{An Example.} We will discuss a phenomenon that is perhaps best introduced with an example: we take the Thomassen graph on 94 vertices (\cite{th} and Fig. 1) and consider the Laplacian $L = D -A$ with eigenvalues ordered as
$$ \lambda_1 \geq \lambda_2 \geq \dots \geq \lambda_{93} >  \lambda_{94} = 0.$$
The graph is 3-regular, the three largest eigenvalues are distinct. Fig. 1 shows the signs of $\phi_2, \phi_3$ (left and middle) and the sign of $\phi_2 \cdot \phi_3$ (right). 

\vspace{-10pt}
\begin{center}
\begin{figure}[h!]
\begin{tikzpicture}[rotate=90]
\node at (0,0) {\includegraphics[width=0.6\textwidth,angle=90]{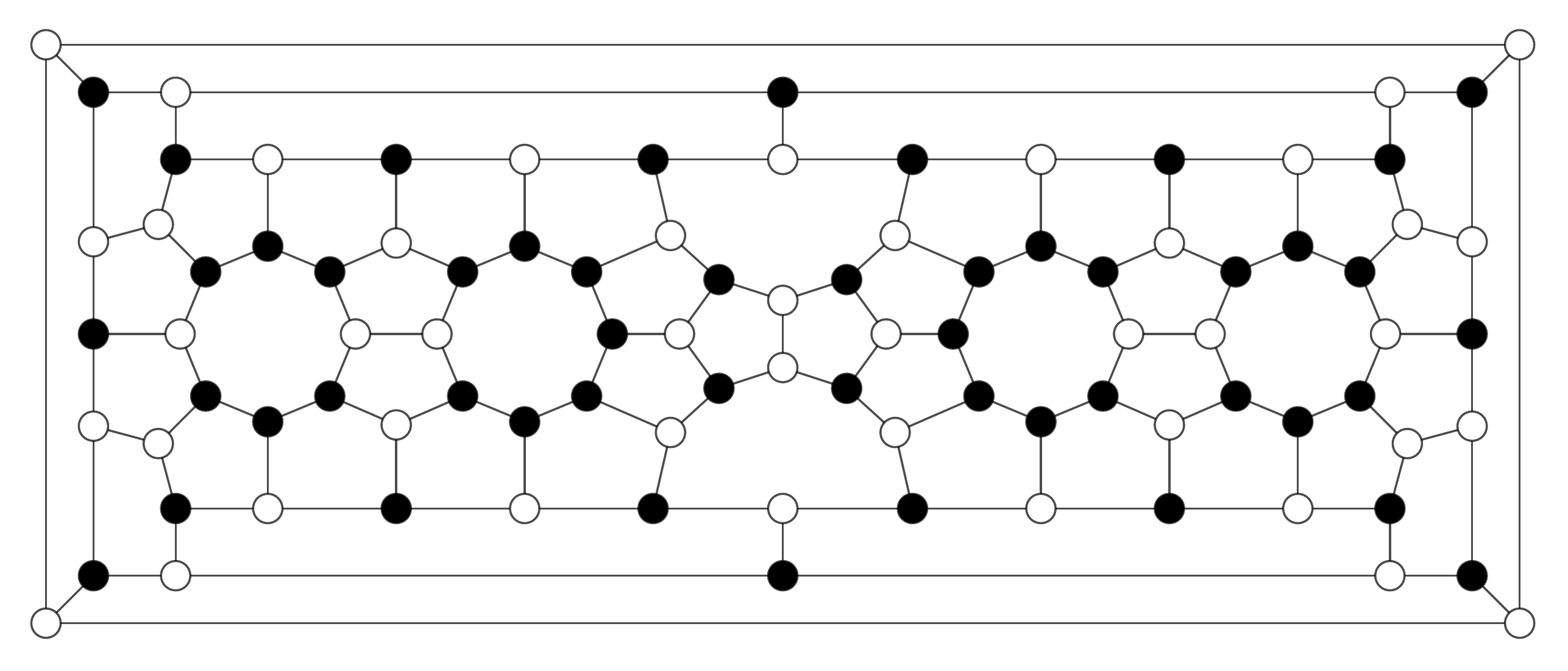}};
\node at (0,-4) {\includegraphics[width=0.6\textwidth,angle=90]{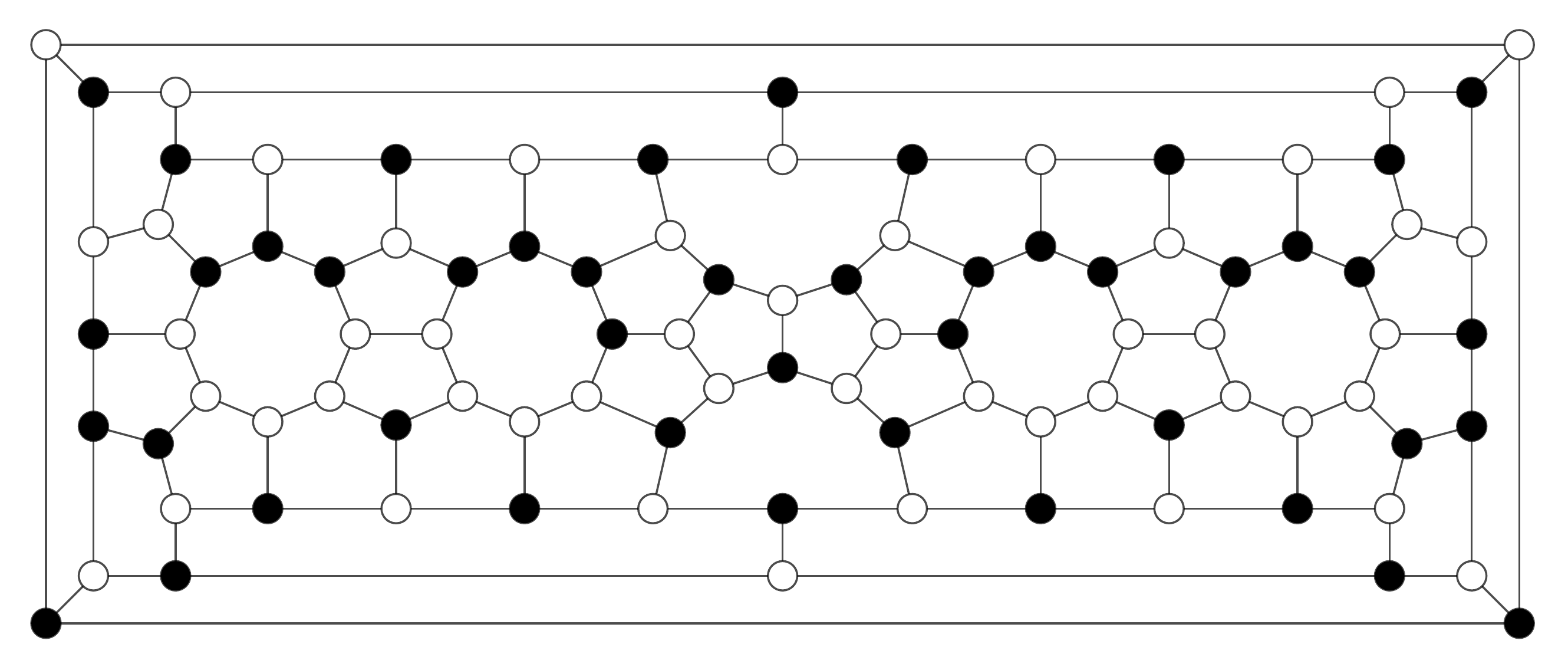}};
\node at (0,-8) {\includegraphics[width=0.6\textwidth,angle=90]{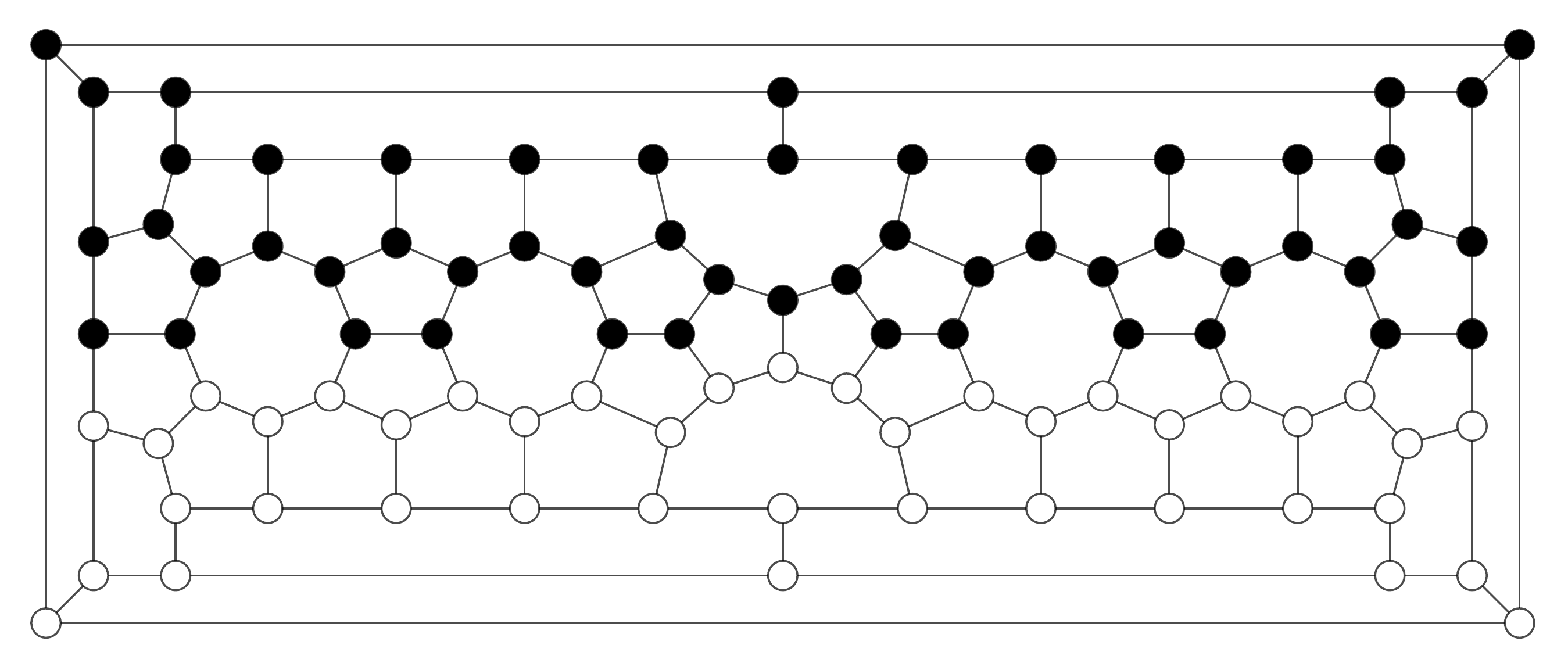}};
\end{tikzpicture}
\vspace{-20pt}
\caption{Sign of the $2^{nd}$ and the $3^{rd}$ eigenvector and of their product.}
\label{fig:1}
\end{figure}
\end{center}

We note that both the second and the third eigenvector have sign changes across most edges: they oscillate essentially as rapidly as the graph allows. In contrast, the (pointwise) product of these high-frequency eigenvectors appears to be much smoother and exhibits a sign pattern typical of low-frequency eigenvectors: positive and negative entries are clustered together and meet across a smooth interface. 
 This can be made quantitative. Having defined the Graph Laplacian as $L = D - A$, we have for all $f:V \rightarrow \mathbb{R}$ a natural measure of `smoothness' of a function
$$ \left\langle f, L f \right\rangle = \sum_{(i,j) \in E} (f(i) - f(j))^2.$$
This is the discrete analogue of the Dirichlet energy $\left\langle f, (-\Delta) f \right\rangle = \int |\nabla f|^2$.
The quadratic form of an eigenfunction is simply its eigenvalue. Here, the quadratic form of the (pointwise) product of the two eigenfunctions is much smaller
$$
  \left\langle \frac{\phi_2}{\| \phi_2\|}, L \frac{\phi_2}{\| \phi_2\|} \right\rangle \sim 5.5 \sim   \left\langle \frac{\phi_3}{\| \phi_3\|}, L \frac{\phi_3}{\| \phi_3\|} \right\rangle  ~~ \mbox{and} ~~ \left\langle \frac{ \phi_2 \cdot \phi_3}{ \| \phi_2 \cdot \phi_3\|}, L~  \frac{\phi_2 \cdot \phi_3}{ \| \phi_2 \cdot \phi_3\|} \right\rangle  \sim 0.5.$$
Note that $\phi_2$ and $\phi_3$ are orthogonal, therefore $\phi_2 \phi_3$ has mean value 0.
At this point one could wonder whether this is simply a coincidence. We illustrate the quadratic form of all pairwise products (all normalized in $\ell^2$) in Fig. 2. We observe three distinct regions: the product of `smooth' (low-frequency) eigenvectors is smooth, the product of a smooth and an oscillatory eigenvector remains oscillatory while the product of two oscillatory eigenvectors tends to be smooth.
\vspace{-10pt}

\begin{center}
\begin{figure}[h!]
\begin{tikzpicture}[]
\node at (0,0) {\includegraphics[width=0.45\textwidth]{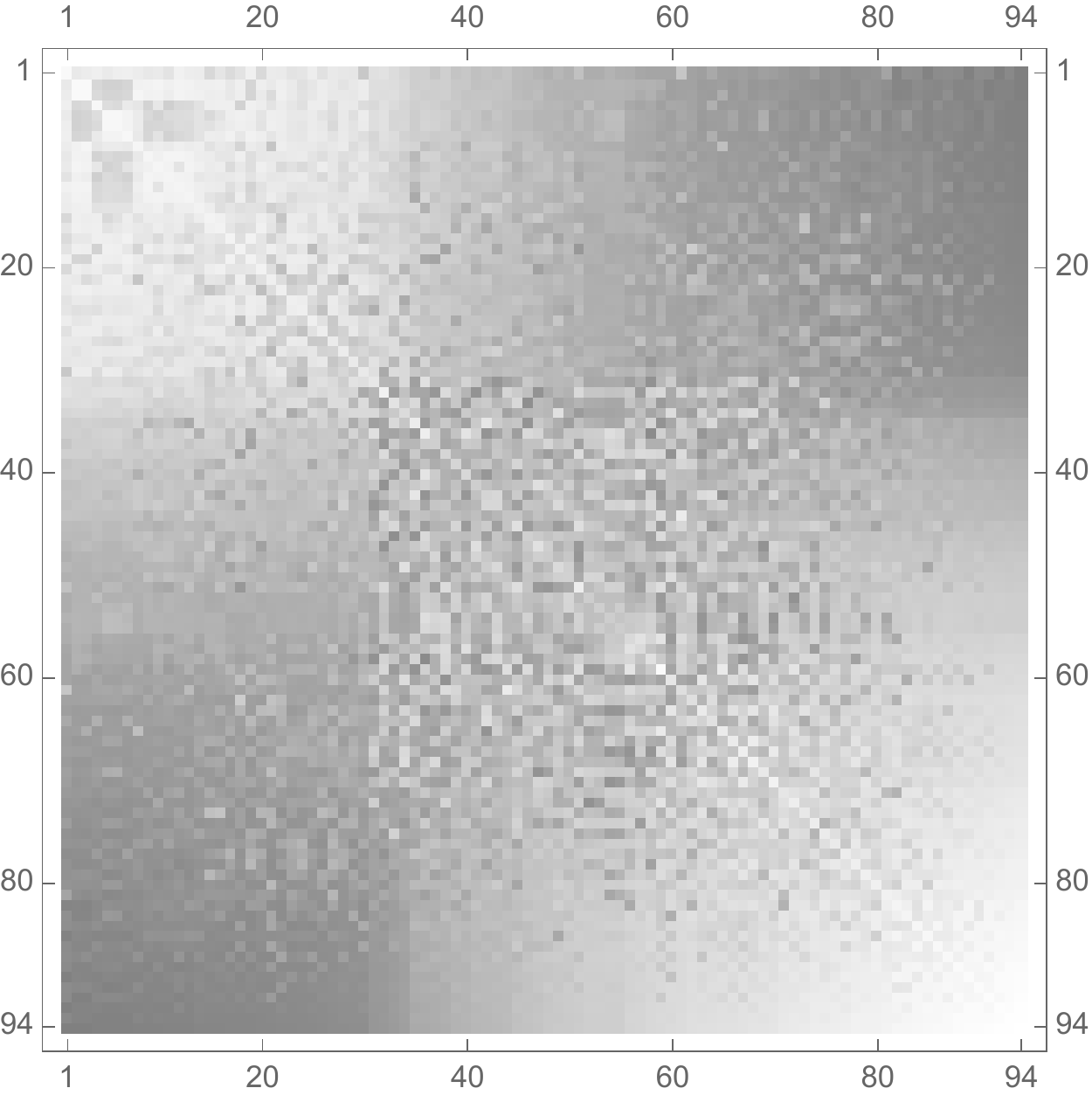}};
\node at (-5, 1) {oscillatory $\cdot$  oscillatory};
\node at (-5, 0.6) {$=$ smooth};
\draw [->] (-4.5, 1.5) -- (-3.5, 2);
\node at (5, 1) {oscillatory $\cdot$  smooth};
\node at (5, 0.6) {$=$ oscillatory};
\draw [->] (4.5, 1.5) -- (3.5, 2);
\node at (-5, -1) {smooth $\cdot$  oscillatory};
\node at (-5, -1.4) {$=$ oscillatory};
\draw [->] (-4.5, -1.6) -- (-3.5, -2);
\node at (5, -1) {smooth $\cdot$  smooth};
\node at (5, -1.4) {$=$ smooth};
\draw [->] (4.5, -1.6) -- (3.5, -2);
\end{tikzpicture}
\caption{Size of the Rayleigh quotient $\left\langle \phi_i \phi_j , L (\phi_i \phi_j)\right\rangle/\| \phi_i \phi_j \|_{\ell^2}^2$ (the lighter, the smaller) on the Thomassen-94 graph \cite{th}.}
\label{fig:2}
\end{figure}
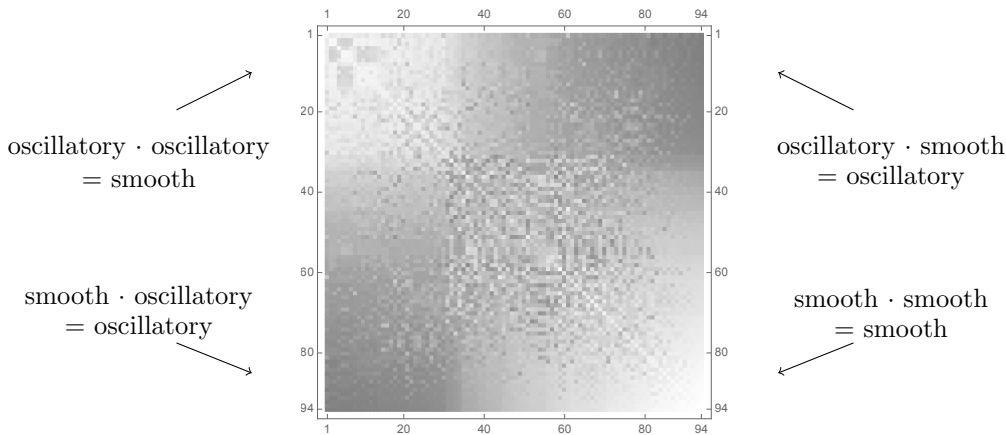
\end{center}
\vspace{-15pt}
The product of two smooth functions is, maybe unsurprisingly, smooth. 
The region $\mbox{smooth} \cdot \mbox{oscillatory} = \mbox{oscillatory}$ is also expected: multiplying a highly oscillatory function with a function varying smoothly across the graph, one can think of the product as a slow modulation of a highly oscillatory function which remains highly oscillatory. We will investigate the fact that the product of two highly oscillatory functions becomes smooth.

\begin{center}
\begin{figure}[h!]
\begin{tikzpicture}[rotate=90]
\node at (0,0) {\includegraphics[width=0.3\textwidth]{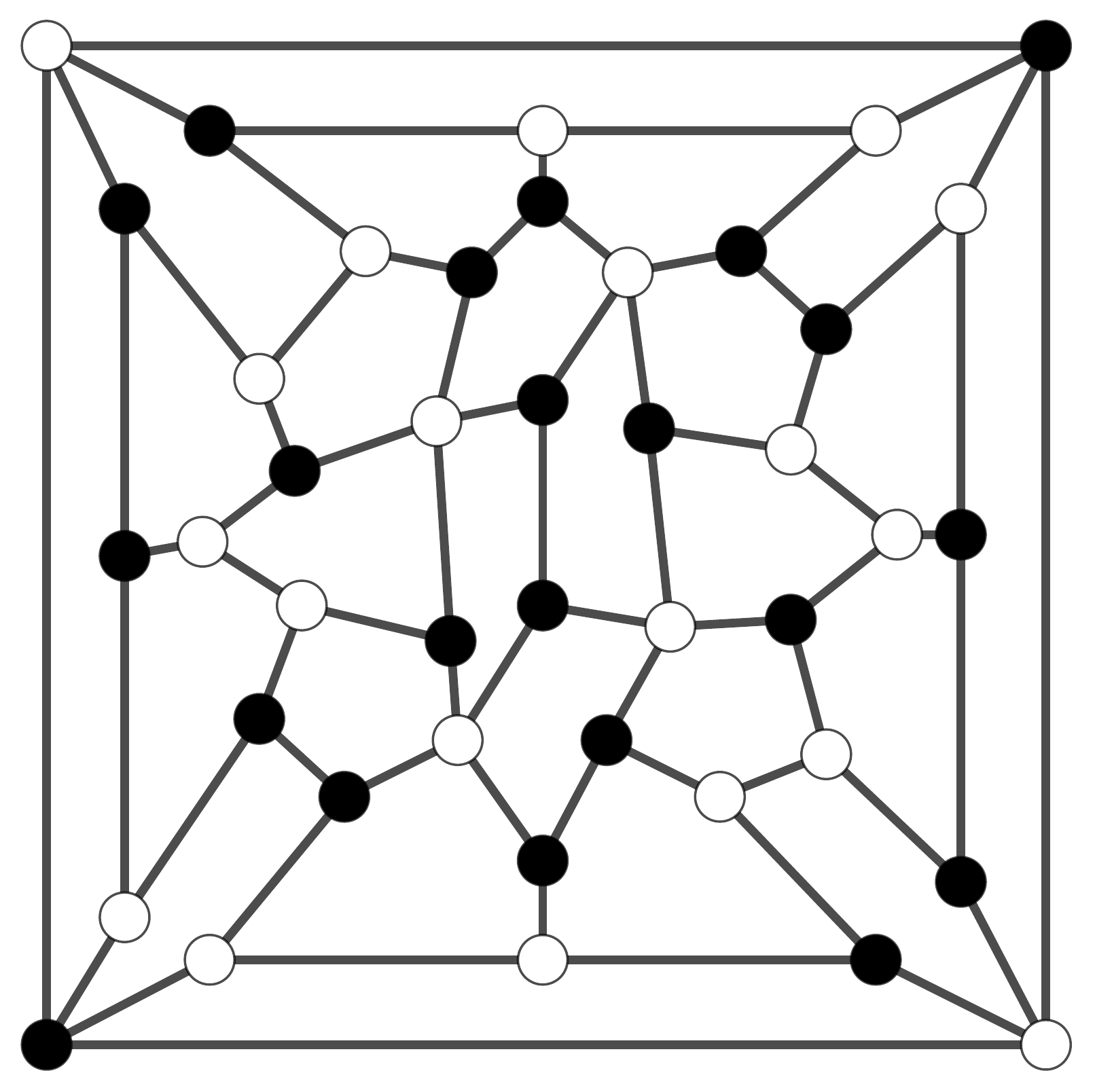}};
\node at (0,-4) {\includegraphics[width=0.3\textwidth]{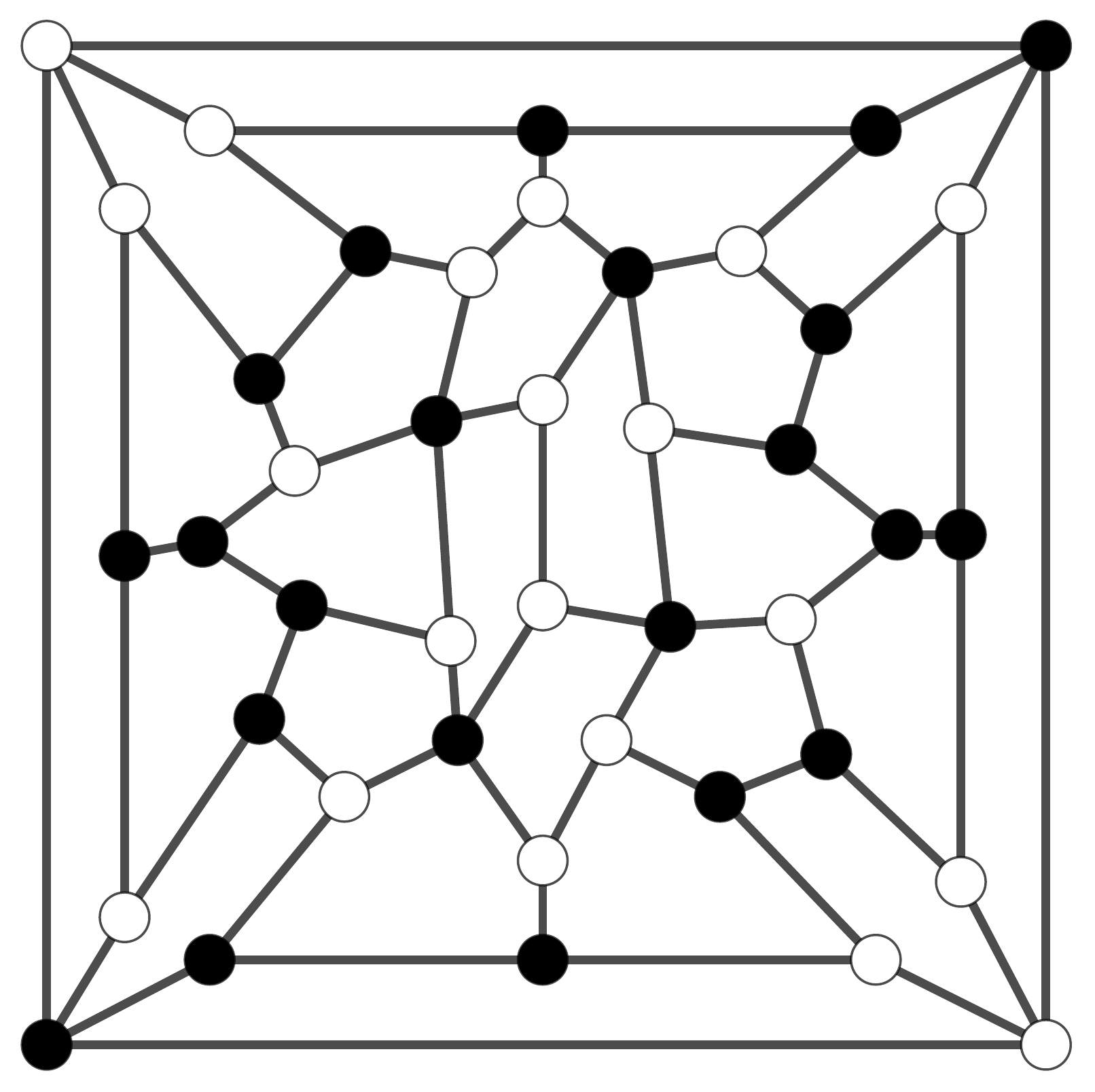}};
\node at (0,-8) {\includegraphics[width=0.3\textwidth]{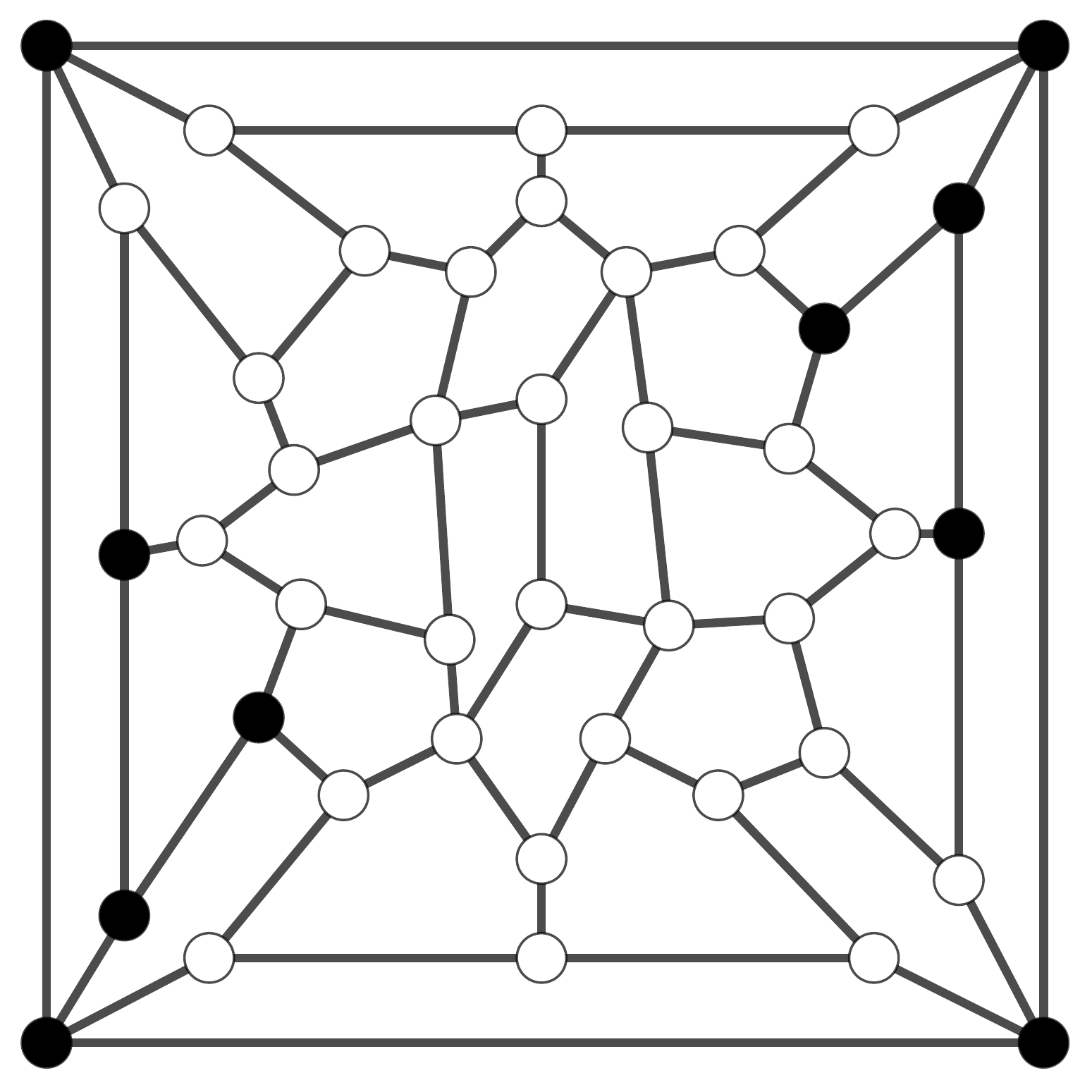}};
\end{tikzpicture}
\vspace{-10pt}
\caption{Sign of the first and second eigenvector on the Wiener-Araya graph \cite{wiener} and the sign of their product.}
\label{fig:3}
\end{figure}
\end{center}

\vspace{-10pt}

\subsection{A Heuristic Explanation.} Let us assume first for simplicity that $G=(V,E)$ is $d-$regular. Then
$$ L = d \cdot \mbox{Id}_{n \times n} - A$$
and the Gerschgorin theorem immediately tells us that $\sigma(L) \subset [0, 2d]$. If $\phi$ is an eigenvector
whose eigenvalue is close to $2d$, then $A \phi \sim - d \cdot \phi$ which means
$$ \frac{1}{d}\sum_{j \sim i} \phi(j) \sim -  \phi(i).$$
The typical value at an adjacent vertex $j$ is close to the negative value in the vertex $i$ itself. Moreover, if $A \phi = -d \cdot \phi$, then $G=(V,E)$ is bipartite and the eigenvector is constant on each component (see Theorem 2 for a stability version of this statement).
The main point of this note is to point out that there is an approximate version for eigenvectors whose eigenvalue is close to $2d$ (it is easy to see that some restriction of this type is necessary: sign cancellation phenomena of this type do not happen on Erd\H{o}s-Renyi random graphs). One way of interpreting this is that highly
oscillatory functions have to `often' change sign across edges, Moreover, the eigenvalue being close to $2d$ forces the graph to have a large bipartite component, so it becomes possible to have sign changes across most edges leading to a consistent pattern (something that would not be possible in the presence of many triangles, for example). The pointwise product of two such functions then leads to cancellation of these oscillatory patterns and results in a `smooth' function.

\begin{center}
\begin{figure}[h!]
\begin{tikzpicture}[rotate=90]
\node at (0,0) {\includegraphics[width=0.31\textwidth]{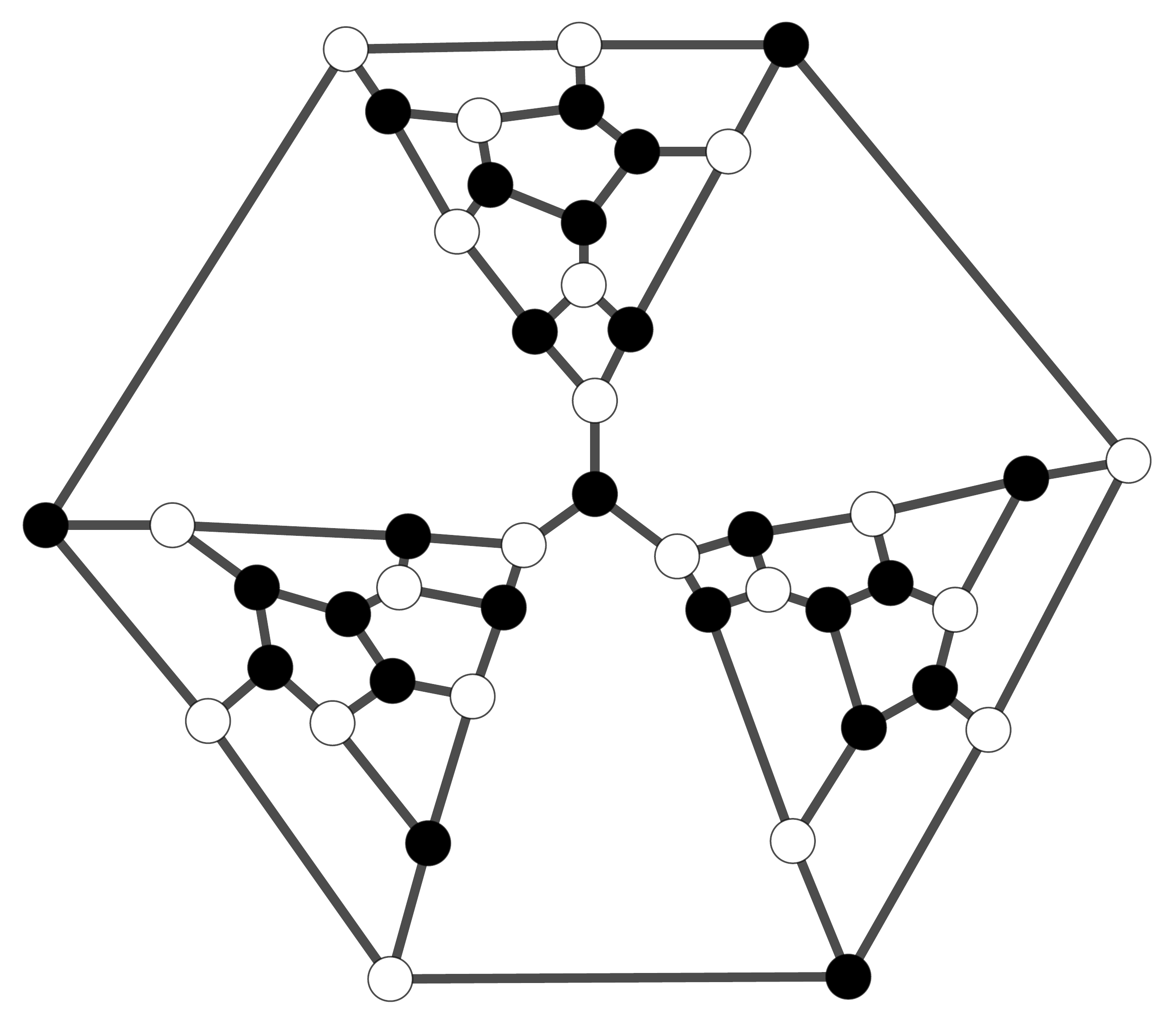}};
\node at (0,-4) {\includegraphics[width=0.31\textwidth]{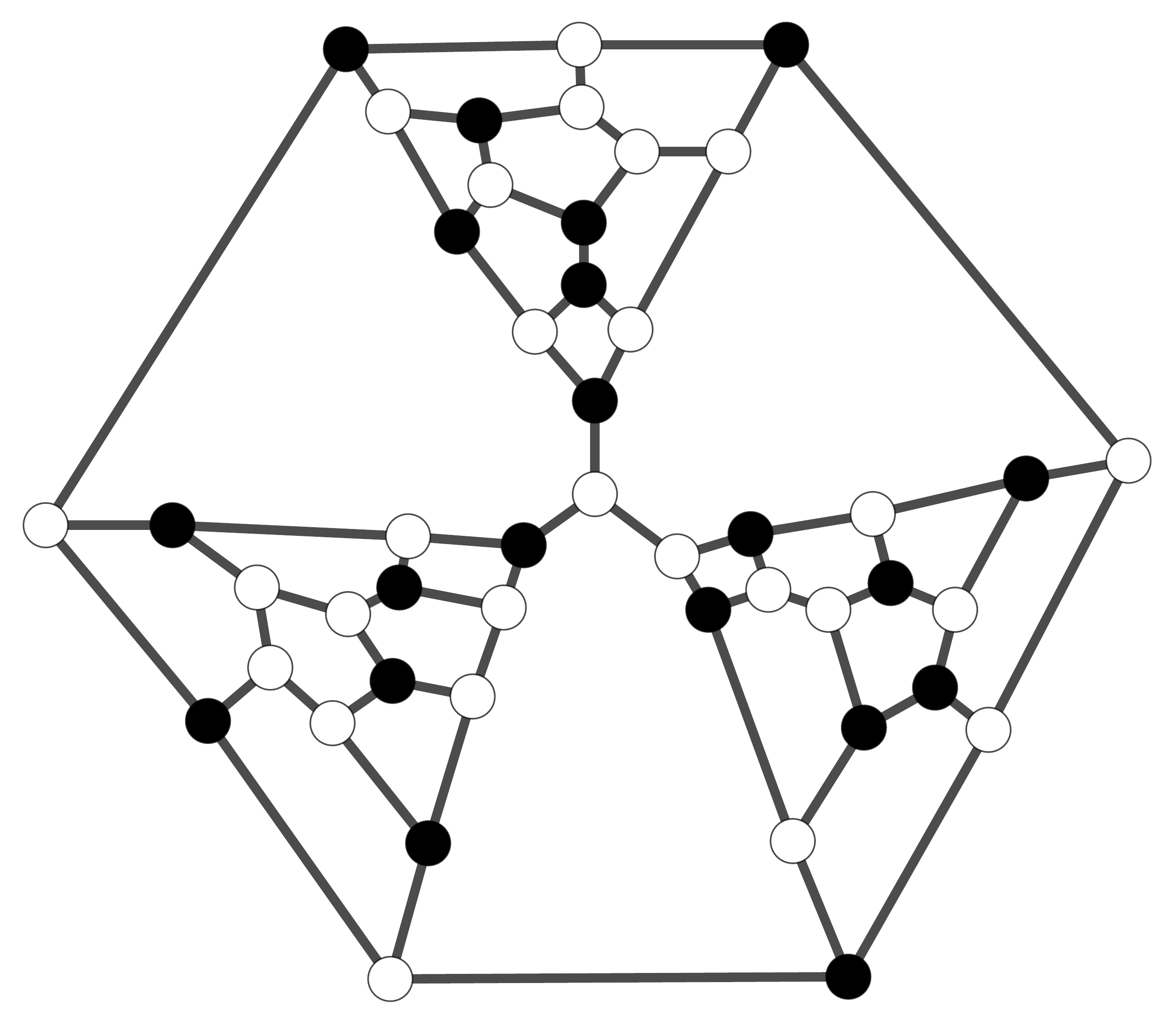}};
\node at (0,-8) {\includegraphics[width=0.31\textwidth]{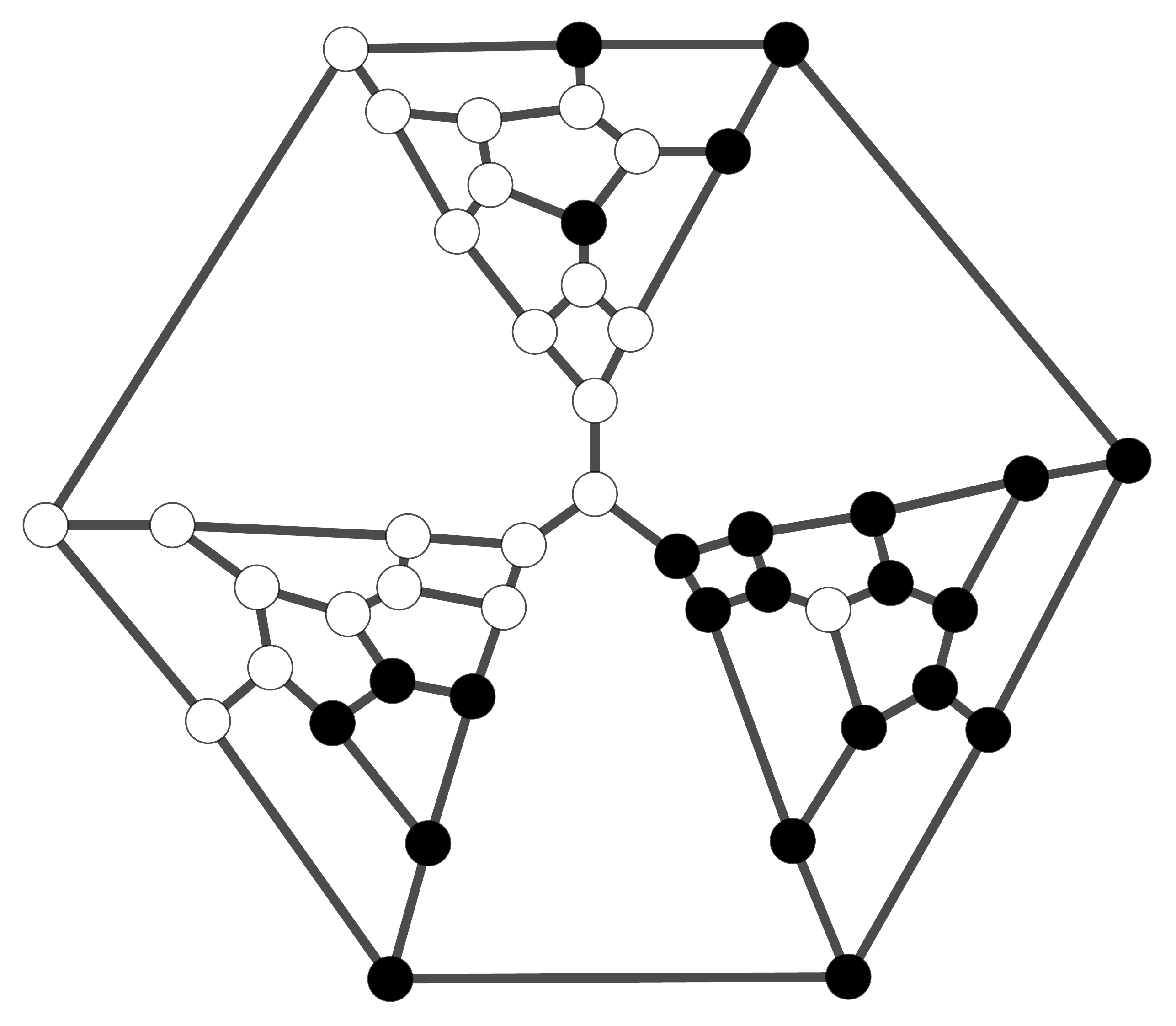}};
\end{tikzpicture}
\vspace{-10pt}
\caption{Sign of the first and second eigenvector on the Tutte graph \cite{tutte} as well as sign of their product.}
\label{fig:4}
\end{figure}
\end{center}

\subsection{Related results.} The question is a priori not meaningful in the continuous setting since there is no such thing as `the most oscillatory function'.  In the continuous setting, the naturally related question is as follows: given a compact manifold/domain with suitable boundary conditions, we obtain a sequence $(\phi_n)_{n=1}^{\infty}$ of $L^2-$normalized Laplacian eigenfunctions. Suppose now that $-\Delta \phi_{\mu} = \mu \cdot \phi_{\mu}$ and $-\Delta \phi_{\lambda} = \lambda \cdot \phi_{\lambda}$, what can be said about the spectral resolution
$$ \phi_{\mu} \cdot \phi_{\lambda} = \sum_{n=1}^{\infty} \left\langle \phi_{\mu} \phi_{\lambda}, \phi_n \right\rangle \phi_n?$$
The term $ \left\langle \phi_{\mu} \phi_{\lambda}, \phi_n \right\rangle$ is also sometimes known as a triple product. It is easy to analyze if the underlying manifold is a torus $\mathbb{T}^d$ but already on the sphere $\mathbb{S}^d$, this is somewhat involved (Clebsch-Gordon coefficients).
Except for the special case where additional structure is present (see  Bernstein \& Reznikoff \cite{bern}, Kr\"otz \& Stanton \cite{krotz}, Sarnak \cite{Sarnak}), there are relatively few results in the literature. The author \cite{stein1} proposed a local interpretation: whether $\phi_{\lambda} \cdot \phi_{\mu}$
significantly shifts in the spectrum depends on whether the local wave-structure is aligned or not. \cite{stein1} mentions the $\mbox{oscillatory} \cdot \mbox{oscillatory} = \mbox{smooth}$ phenomenon (without explanation).

\vspace{-20pt}
\begin{center}
\begin{figure}[h!]
\begin{tikzpicture}[rotate=90]
\node at (0,0.25) {\includegraphics[width=0.31\textwidth]{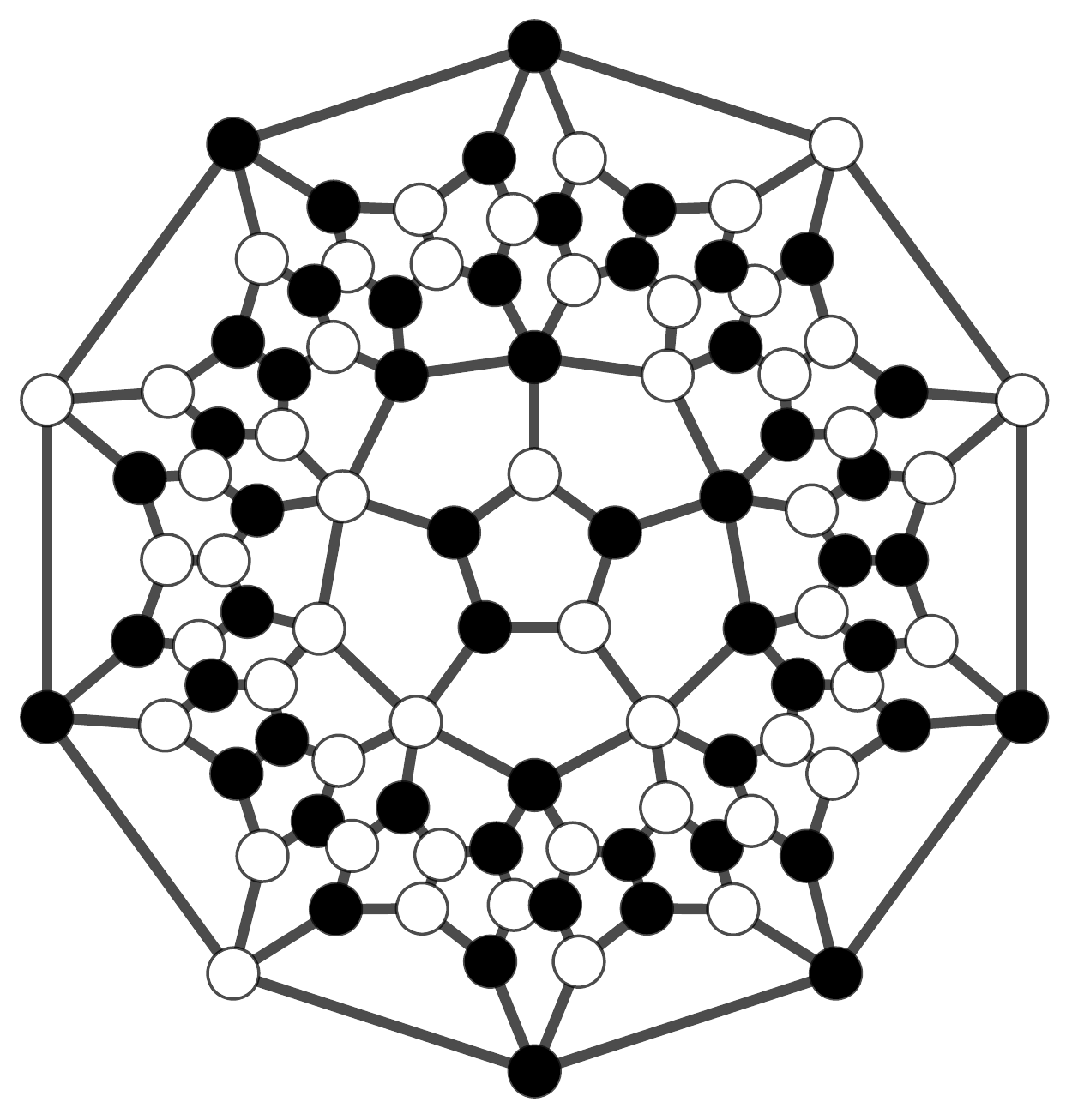}};
\node at (0,-4) {\includegraphics[width=0.31\textwidth]{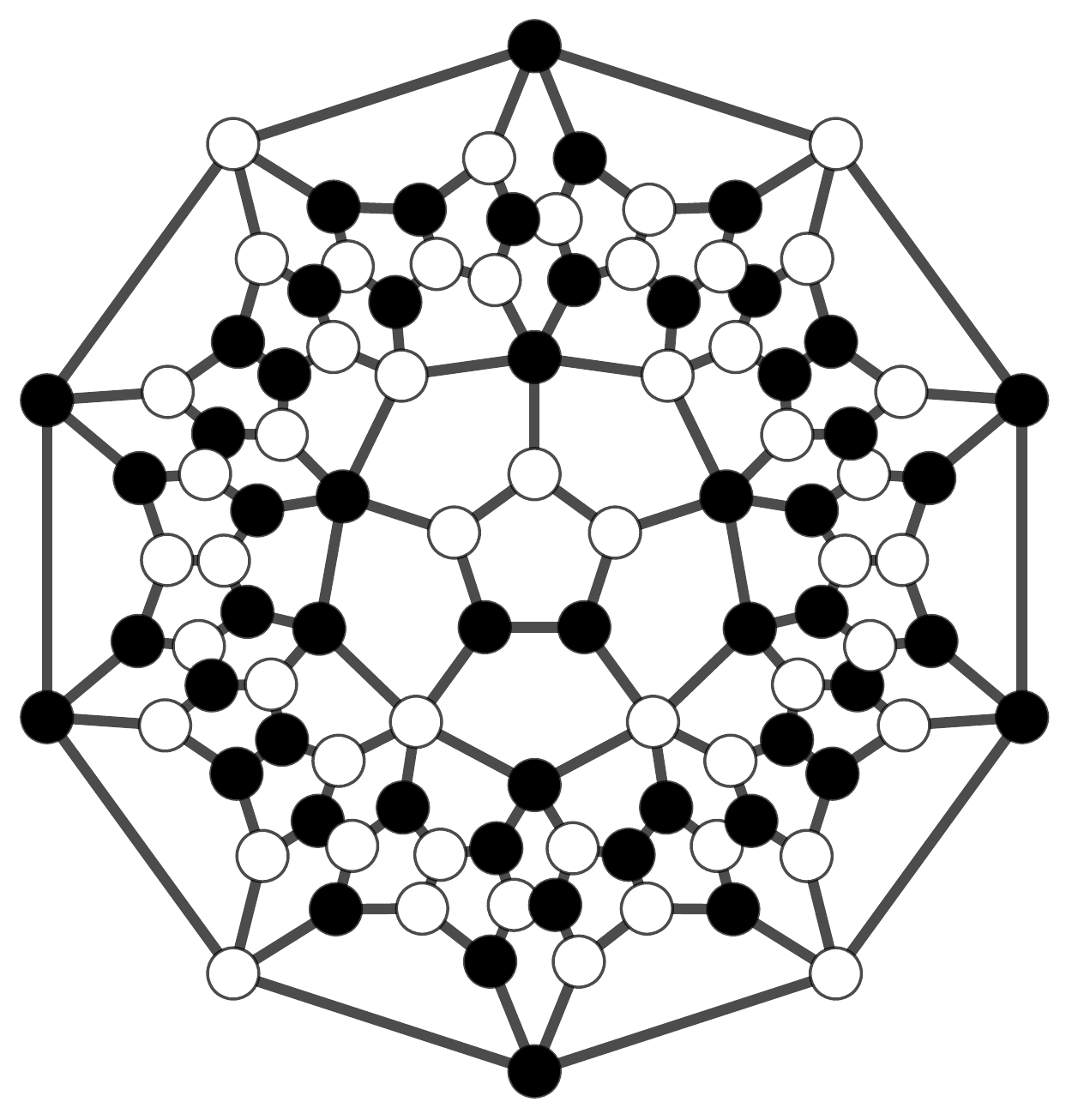}};
\node at (0,-8.25) {\includegraphics[width=0.31\textwidth]{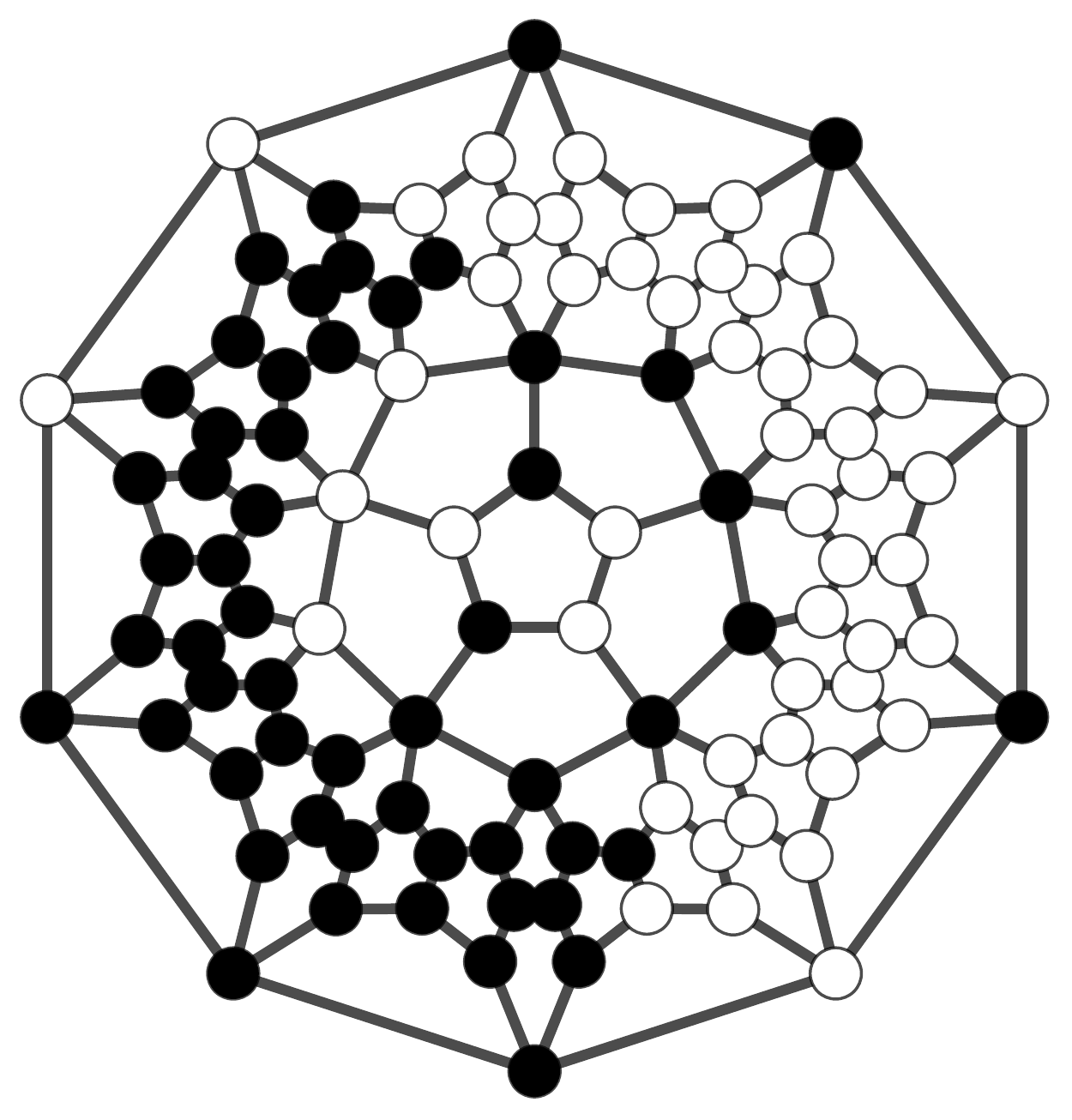}};
\end{tikzpicture}
\vspace{-10pt}
\caption{Sign of the $18^{th}$ and $19^{th}$ eigenvector on Thomassen-105 (see \cite{th0}) and the sign of their product. This example deviates from the others insofar as the two eigenvectors are not `very' extremal -- however, there is nonetheless a cancellation in their sign pattern and the product has a much smaller quadratic form.}
\label{fig:8}
\end{figure}
\end{center}
\vspace{-20pt}

Lu, Sogge and the author \cite{lu} gave a quantitative description of the $\mbox{smooth} \cdot \mbox{smooth} = \mbox{smooth}$ regime in the continuous setting (see also Jin \cite{jin} and Wyman \cite{wyman}). The product $\phi_{\lambda} \cdot \phi_{\mu}$ appear in a variety of different applications such as finding appropriate eigenfunctions for the dimensionality reduction of high-dimensional data (see Cloninger \& Steinerberger \cite{cloninger},
Kohli, Cloninger \& Mishne \cite{kohli}),  in numerical aspects of the Kohn-Sham density functional theory (see Lin, Lu \& Ying \cite{lin}) and in problems related to shape matching (see Litany, Rodola, Bronstein \& Bronstein \cite{bronstein}).
 The largest eigenvectors (and the associated eigenvalues) are known to have great significance in combinatorics, specifically the chromatic number and size of independent set (i.e. the Hoffman bound \cite{hoffman}, the work of Godsil \& Newman \cite{godsil}), the Max-Cut problem (see e.g. Delorme \& Pojak \cite{delorme}, Mohar \& Poljak \cite{mohar}, Trevisan \cite{trevisan}), the domination number (see e.g. Nikiforov \cite{niki}) and many others, however, we are not aware of any results regarding pointwise products of eigenvectors.

\section{Statement of Results} We give a relatively simple quantitative explanation of the result for the Laplacian $L = D - A$. However, it should be emphasized
that the underlying phenomenon is more general than that and should have analogues for other types of Graph Laplacians and may also have other quantitative
formulations for $L=D-A$.
\subsection{A simple bound on the Dirichlet energy} The Gerschgorin bound implies
$$ 0 \leq D-A \leq 2D$$
and we will call an eigenvalue `big' if $(D-A) \phi$ is close to $2D \phi$. Alternatively, we can assume that 
$$ \left[ 2D - (D-A) \right] \phi = (D+A) \phi \qquad \mbox{is small.}$$
There is a simple bound explaining this phenomenon for delocalized functions.
\begin{theorem} Let $G$ be a graph and let $\phi, \psi \in V \rightarrow \mathbb{R}$ be two functions satisfying
$$ \left\langle \phi, (D+A) \phi \right\rangle \leq \varepsilon \| \phi\|^2_{\ell^2} \quad \mbox{and} \quad  \left\langle \psi, (D+A) \psi \right\rangle \leq \varepsilon \| \psi\|^2_{\ell^2}.$$
Then their pointwise product $\phi \cdot \psi: V \rightarrow \mathbb{R}$ satisfies
$$ \left\langle \frac{\phi \cdot \psi}{\| \phi \cdot \psi \|_{\ell^2}}, L \frac{ \phi \cdot \psi}{\| \phi \cdot \psi \|_{\ell^2}} \right\rangle \leq  2\varepsilon \left( \frac{\|\phi\|_{\ell^{\infty}}^2   \|\psi\|_{\ell^2}^2}{\| \phi \cdot \psi \|_{\ell^2}} +  \frac{    \|\phi\|_{\ell^2}^2 \|\psi\|_{\ell^{\infty}}^2 }{\| \phi \cdot \psi \|_{\ell^2}}\right).$$
\end{theorem}

We do not need to assume that $\phi, \psi$ are eigenfunctions or even that they are orthogonal, the choice $\phi = \psi$ is allowed. The cycle graph $C_n$ shows that the inequality is sharp up to at most a factor of 2.
 The most important question when trying to understand the quality of the estimate is how much is lost in the application of H\"older's inequalities
$ \| \phi \cdot \psi \|_{\ell^2} \leq \| \phi \|_{\ell^{\infty}} \cdot \| \psi \|_{\ell^2}$ and $ \| \phi \cdot \psi \|_{\ell^2} \leq \| \psi \|_{\ell^{\infty}} \cdot \| \phi \|_{\ell^2}$. Inequalities of this type are very lossy in general, however, this is not the case here: extremal eigenvectors tend to be `flat' in the sense that their maximal entry is not too much larger than a typical entry. A reason for the tendency towards delocalization is given in \S 2.2 and \S 2.3. It is an interesting question whether other quantitative formulations are possible (both for $L = D-A$ and other Laplacians).
\vspace{0pt}
\begin{center}
\begin{figure}[h!]
\begin{tikzpicture}[rotate=90]
\node at (0,0) {\includegraphics[width=0.31\textwidth]{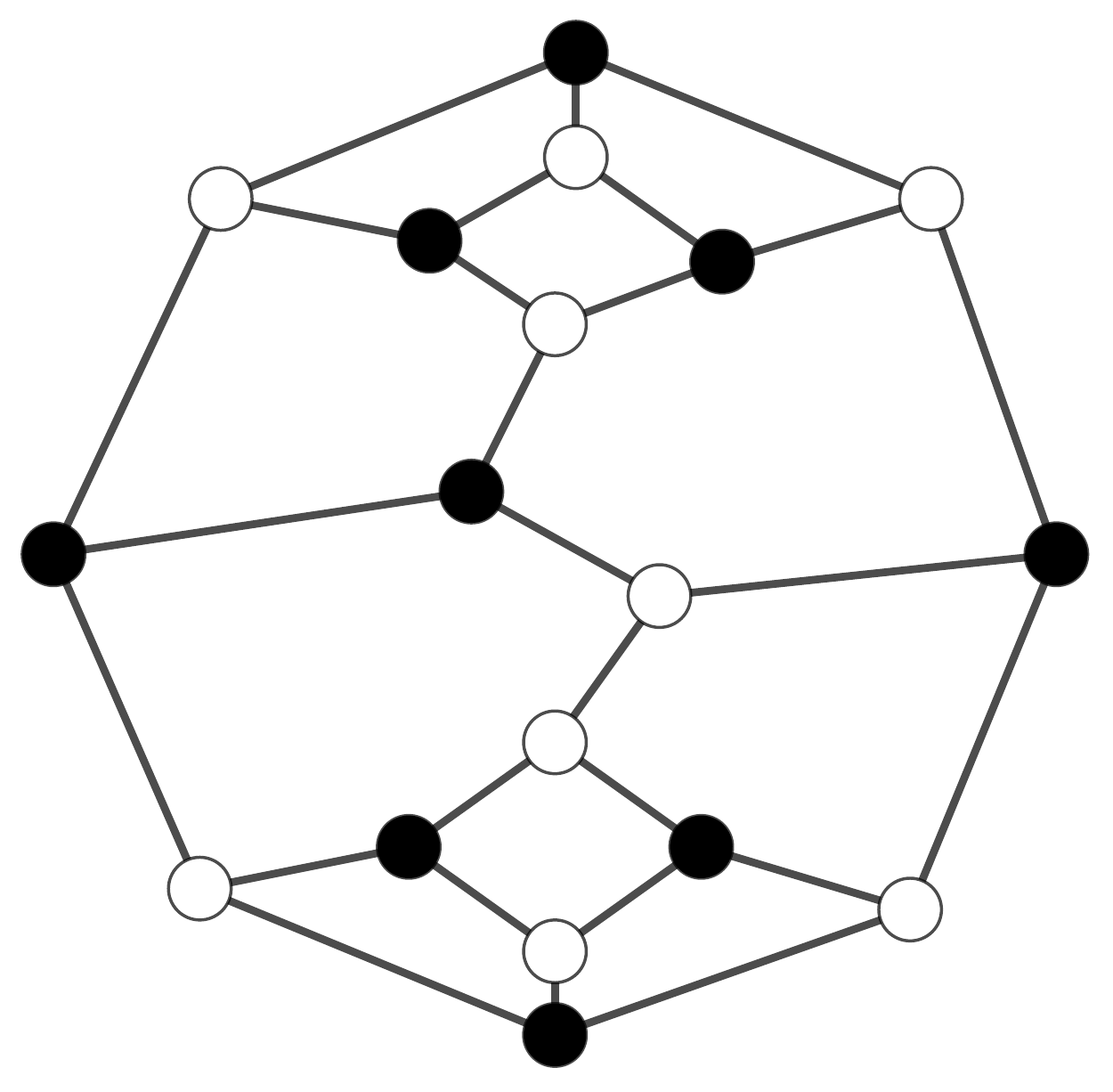}};
\node at (0,-4) {\includegraphics[width=0.31\textwidth]{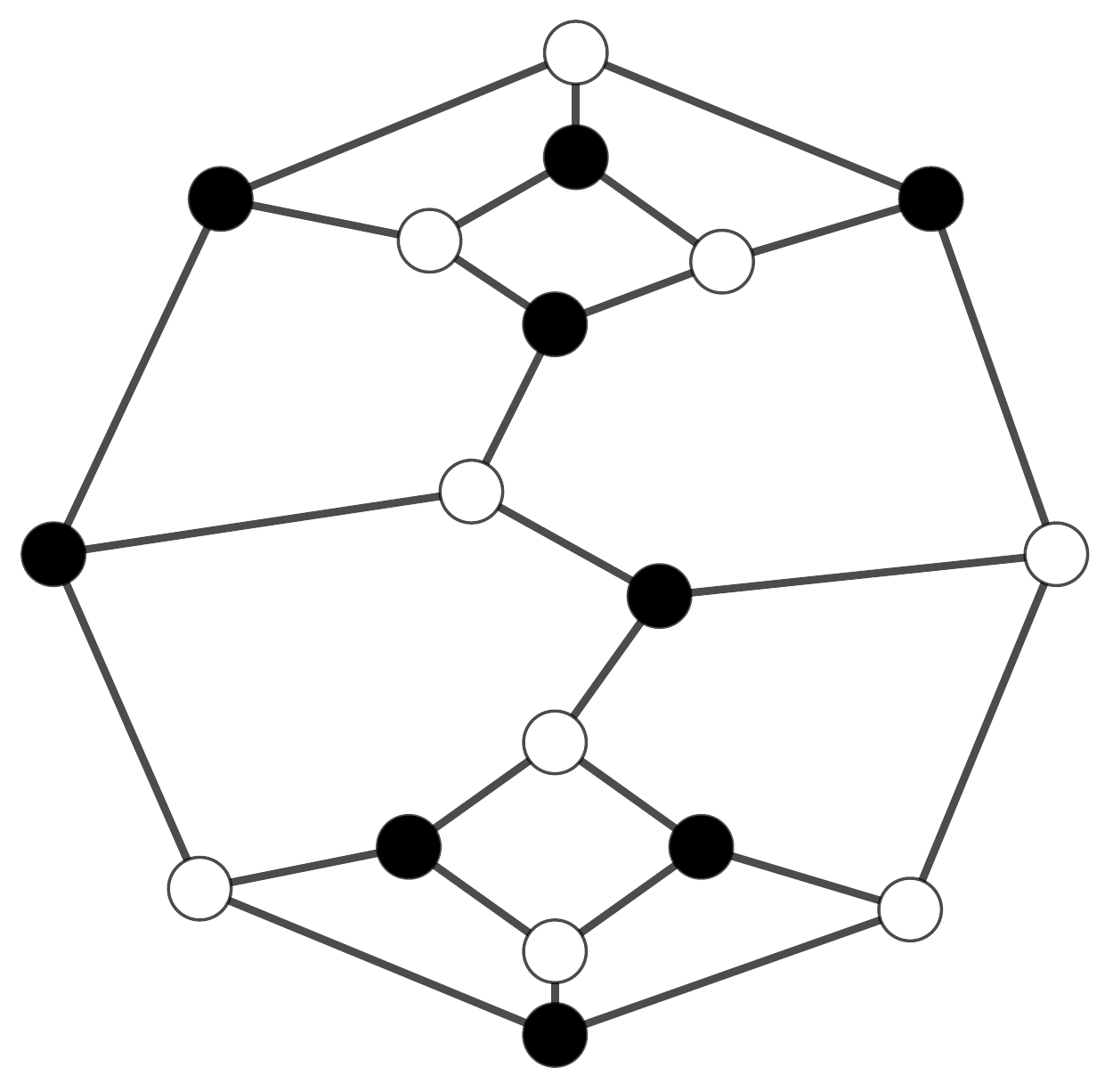}};
\node at (0,-8) {\includegraphics[width=0.31\textwidth]{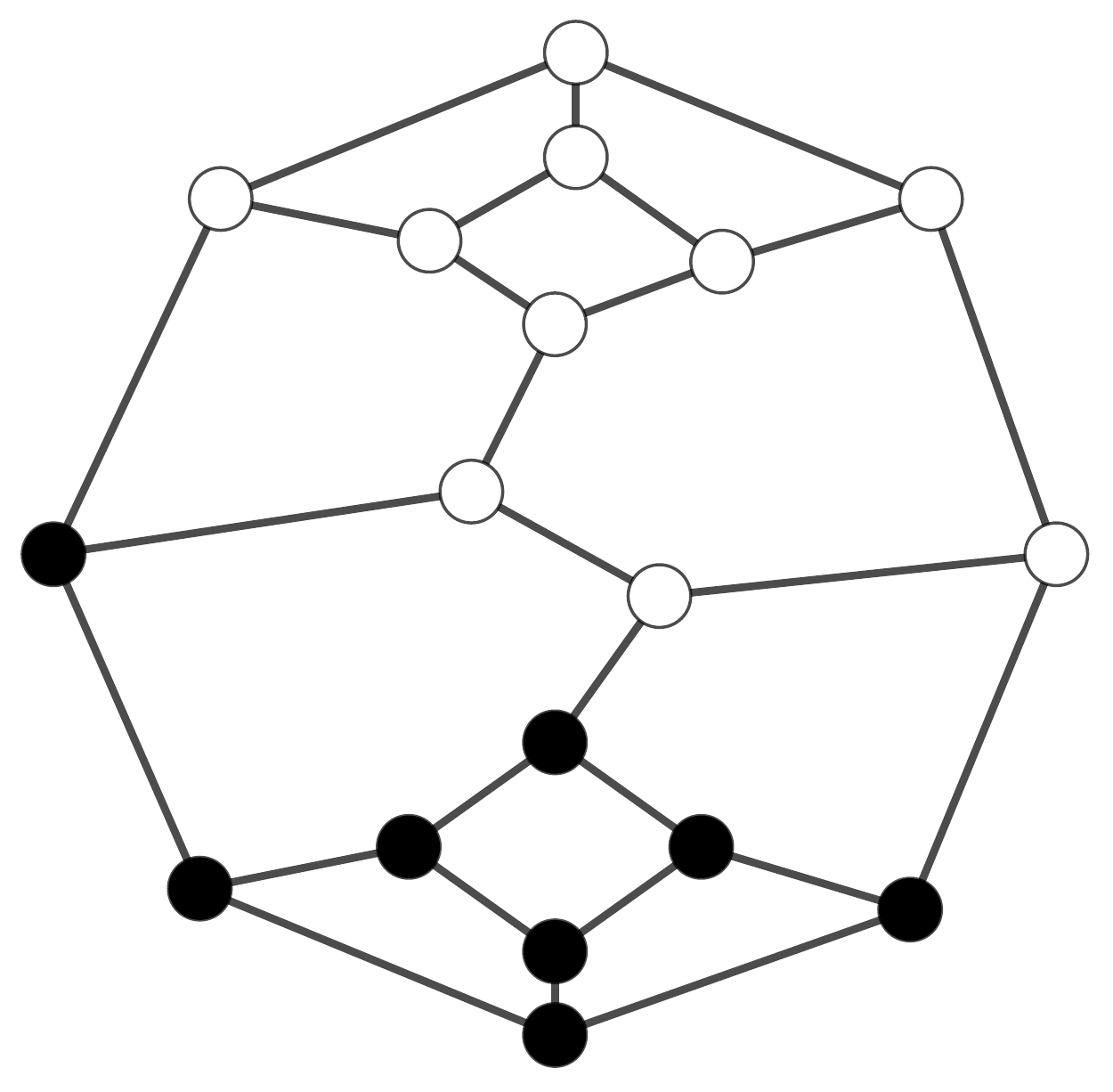}};
\end{tikzpicture}
\vspace{-10pt}
\caption{Sign of the first and second eigenvector on a cubic graph and the sign of their product.}
\label{fig:6}
\end{figure}
\end{center}
\vspace{0pt}

\subsection{A bound on the $\ell^{\infty}-$norm} We prove that for regular graphs the $\ell^{\infty}-$norm of an eigenfunction corresponding to a large eigenvalue, $(D+A)\phi = \varepsilon \phi$, is `small' in a suitable sense. We recall the  matrix $AD^{-1}$: if $\psi:V \rightarrow \mathbb{R}$ is a probability distribution over $V$, then $AD^{-1} \psi$ is the probability distribution after starting with $\psi$ and jumping to a randomly chosen adjacent vertex.\\

We prove a bound stating that if the graph has the property that the random walk started in an arbitrary vertex is not particularly likely to be in any other vertex after a certain number of steps, then this forces the $\ell^{\infty}-$norm to be small.

\begin{theorem} Let $G$ be a $d-$regular graph and let $\phi:V \rightarrow \mathbb{R}$ satisfy
$$ (d \cdot \emph{Id}_{n \times n}+A)\phi = \varepsilon \phi.$$
Suppose $\phi(m) = \max_{v \in V}|\phi(v)|$. Then, for any $k \in \mathbb{N}_{\geq 0}$
$$ \| \phi\|_{\ell^{\infty}} \leq \left( \frac{d}{d-\varepsilon} \right)^{2k}  \left(\sum_{j \in V} \left[(AD^{-1})^{2k} \delta_m\right](j)^2\right)^{1/2}\|\phi\|_{\ell^2}. $$
\end{theorem}

We can illustrate the result in a toy example. Suppose $G$ is a $d-$regular graph on $2n$ vertices that is `mostly' bipartite: there is a decomposition $V = A \cup B$ with $|A| = n = |B|$ such that only $k$ edges do not run between $A$ and $B$.  Using the Rayleigh-Ritz formulation, we see that there exists a vector $\phi$ satisfying
$$ (d \cdot \mbox{Id}_{n \times n}+A)\phi = \varepsilon \phi \qquad \mbox{with} \qquad 0 \leq \varepsilon \leq \frac{2k}{n}.$$
At the same time, for a generic `random' graph of this type, we expect that a random walk jumps mostly from partition to partition but equidistributes within each partition within $2k \sim \log{n}$ steps. This, for $k$ sufficiently small, results in a bound $\| \phi\|_{\ell^{\infty}} \lesssim n^{-1/2} \|\phi\|_{\ell^2}$ which is clearly best possible.

\begin{center}
\begin{figure}[h!]
\begin{tikzpicture}[rotate=90]
\node at (0,0) {\includegraphics[width=0.24\textwidth]{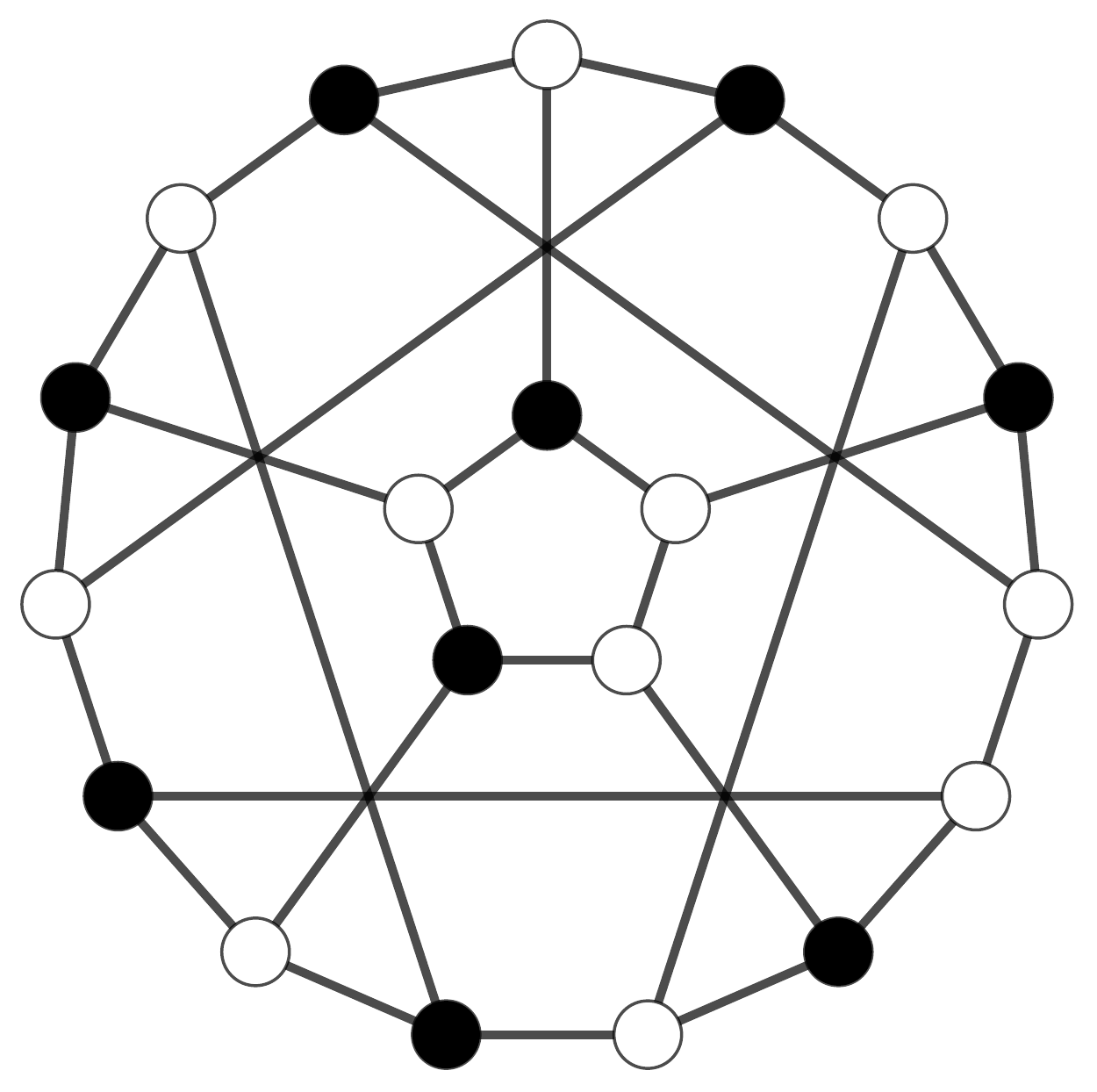}};
\node at (0,-4) {\includegraphics[width=0.24\textwidth]{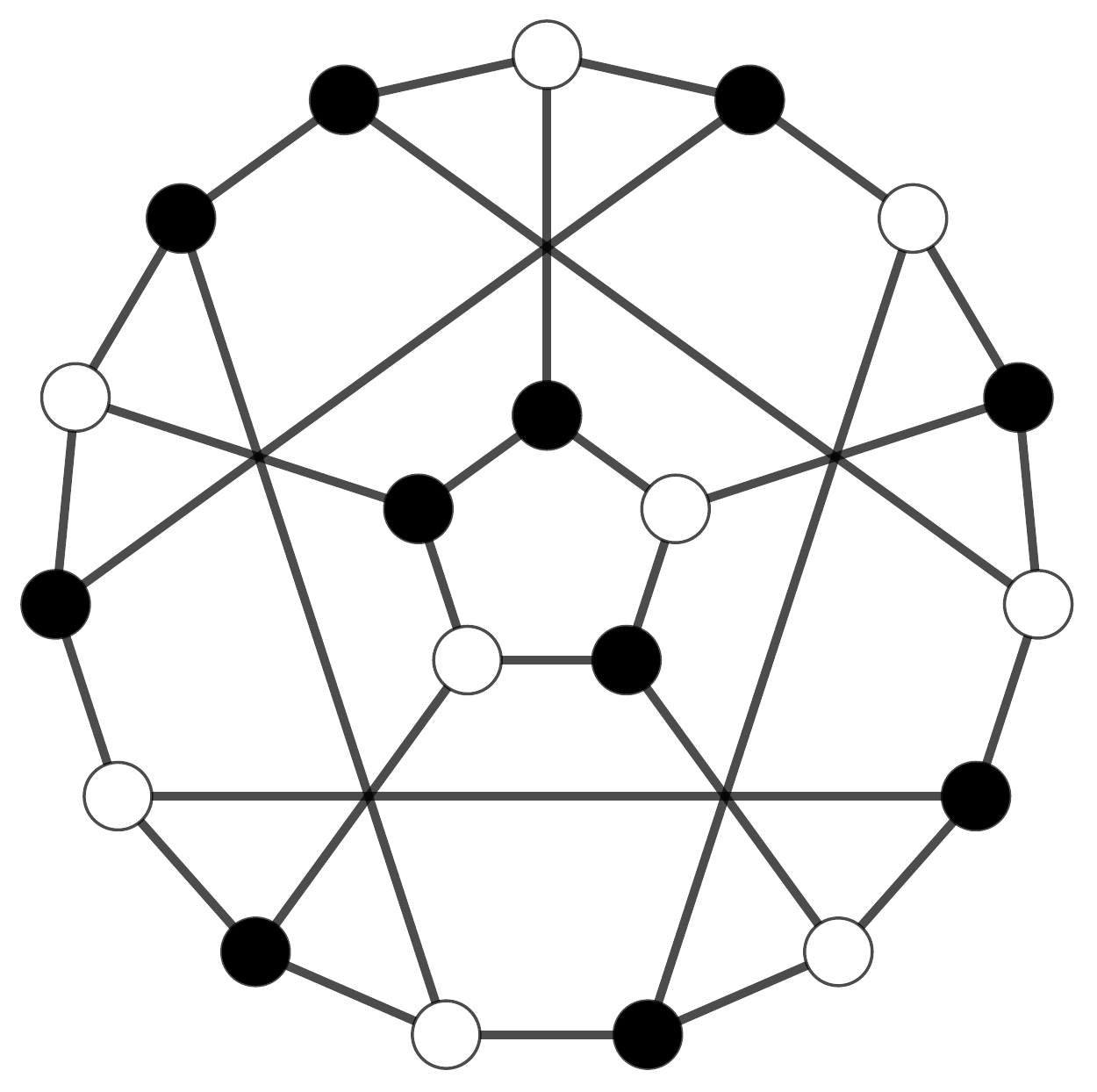}};
\node at (0,-8) {\includegraphics[width=0.24\textwidth]{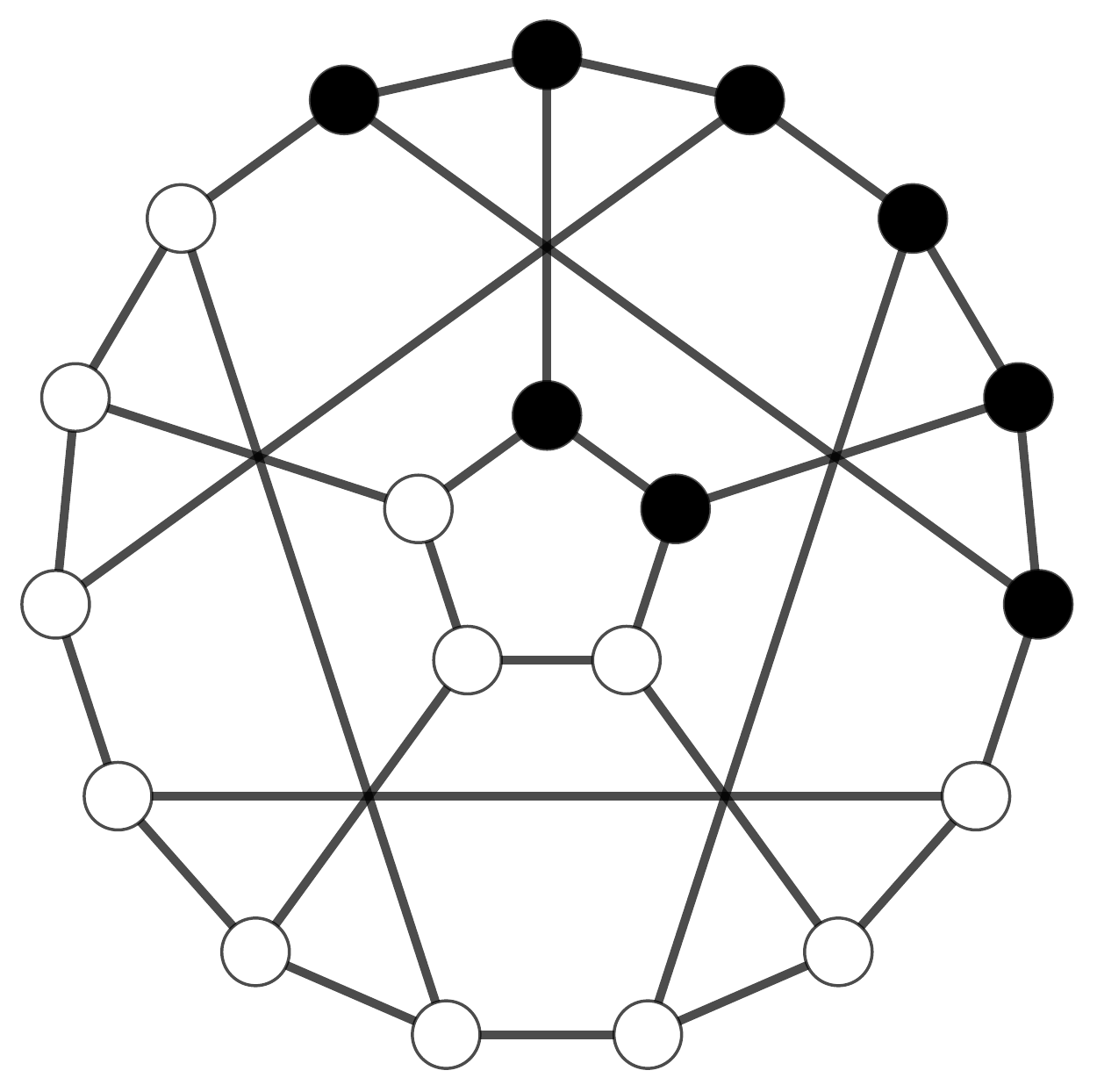}};
\end{tikzpicture}
\vspace{-10pt}
\caption{Sign of the first and second eigenvector on the Flower Snark $J_5$ \cite{isaacs} as well as the sign of their pointwise product.}
\label{fig:7}
\end{figure}
\end{center}

\subsection{Another perspective} We conclude with another elementary estimate. 
\begin{corollary} Let $G$ be a $d-$regular graph and let $\phi:V \rightarrow \mathbb{R}$ satisfy
$$ (D+A)\phi = \varepsilon \phi.$$
Then
 $$ \left\| \phi - (AD^{-1})^2 \phi \right\|_{\ell^2} = \frac{\varepsilon}{d}\left(2 - \frac{\varepsilon}{d}\right)  \| \phi\|_{\ell^2}.$$
\end{corollary}
Recall that $AD^{-1}$ corresponds to a random walk. From $(D+A)\phi = \varepsilon \phi$, we conclude
that $(\mbox{Id}_{n \times n} + AD^{-1}) \phi = (\varepsilon/d) \phi$ and therefore
$$ \left\| \phi - (AD^{-1}) \phi \right\|_{\ell^2} = \left( 2 - \frac{\varepsilon}{d} \right) \left\| \phi \right\|_{\ell^2}.$$
This means that one step of the random walk is very different from the function (almost as different as possible);
however, as indicated by the Corollary, two steps of the random walk barely change anything at all. This leads to
another interpretation: such eigenvectors are smooth with respect to random
walks of size 2 and the large bipartite component of the graph respects this type of smoothness. In that sense, the `oscillatory $\cdot$ oscillatory = smooth' regime may be naturally understood
as an extension of `smooth $\cdot$ smooth = smooth' with respect to step size 2.

\section{Proofs}
\subsection{Proof of Theorem 1}
\begin{proof}
The proof is simple: we have
\begin{align*}
\left\langle \phi \cdot \psi, L(\phi \cdot \psi )\right\rangle &=  \sum_{e \in E} (\phi(i) \psi(i) - \phi(j) \psi(j))^2\\
 &=  \sum_{e \in E} (\phi(i) \psi(i) + \phi(i) \psi(j) - \phi(i) \psi(j) - \phi(j) \psi(j))^2
 \end{align*}
 We use $(a-b)^2 \leq 2(a^2 + b^2)$ to conclude
 \begin{align*}
\left\langle \phi \cdot \psi, L \phi \cdot \psi \right\rangle &\leq 2  \sum_{e \in E} (\phi(i) \psi(i) + \phi(i) \psi(j))^2 + 2\sum_{e \in E}  (\phi(i) \psi(j) + \phi(j) \psi(j))^2 \\
 &\leq 2 \| \phi\|^2_{\ell^{\infty}}  \sum_{e \in E} (\psi(i) + \psi(j))^2 + 2 \| \psi\|^2_{\ell^{\infty}}\sum_{e \in E}  (\phi(i) + \phi(j) )^2
\end{align*}
These sums simplify since
\begin{align*}
\sum_{e \in E} ( \psi(i) + \psi(j))^2  &= \sum_{e \in E} \psi(i)^2 + \psi(j)^2 +  2\psi(i) \psi(j) \\
&= \sum_{i \in V} \deg(i) \psi(i)^2 + \sum_{e \in E} 2 \psi(i) \psi(j) \\
&= \left\langle \psi, D \psi \right\rangle + \left\langle \psi, A \psi \right\rangle \leq \varepsilon \| \psi\|^2_{\ell^2}.
\end{align*}
Arguing in the same way for $\phi$, we arrive at 
\begin{align*}
\left\langle \phi \cdot \psi, L(\phi \cdot \psi )\right\rangle \leq 2\varepsilon \left( \|\phi\|_{\ell^{\infty}}^2  \cdot \|\psi\|_{\ell^2}^2 +   \|\phi\|_{\ell^2}^2 \cdot \|\psi\|_{\ell^{\infty}}^2  \right).
\end{align*}
\end{proof}

\vspace{-25pt}
\begin{center}
\begin{figure}[h!]
\begin{tikzpicture}[rotate=90]
\node at (0,0) {\includegraphics[width=0.31\textwidth]{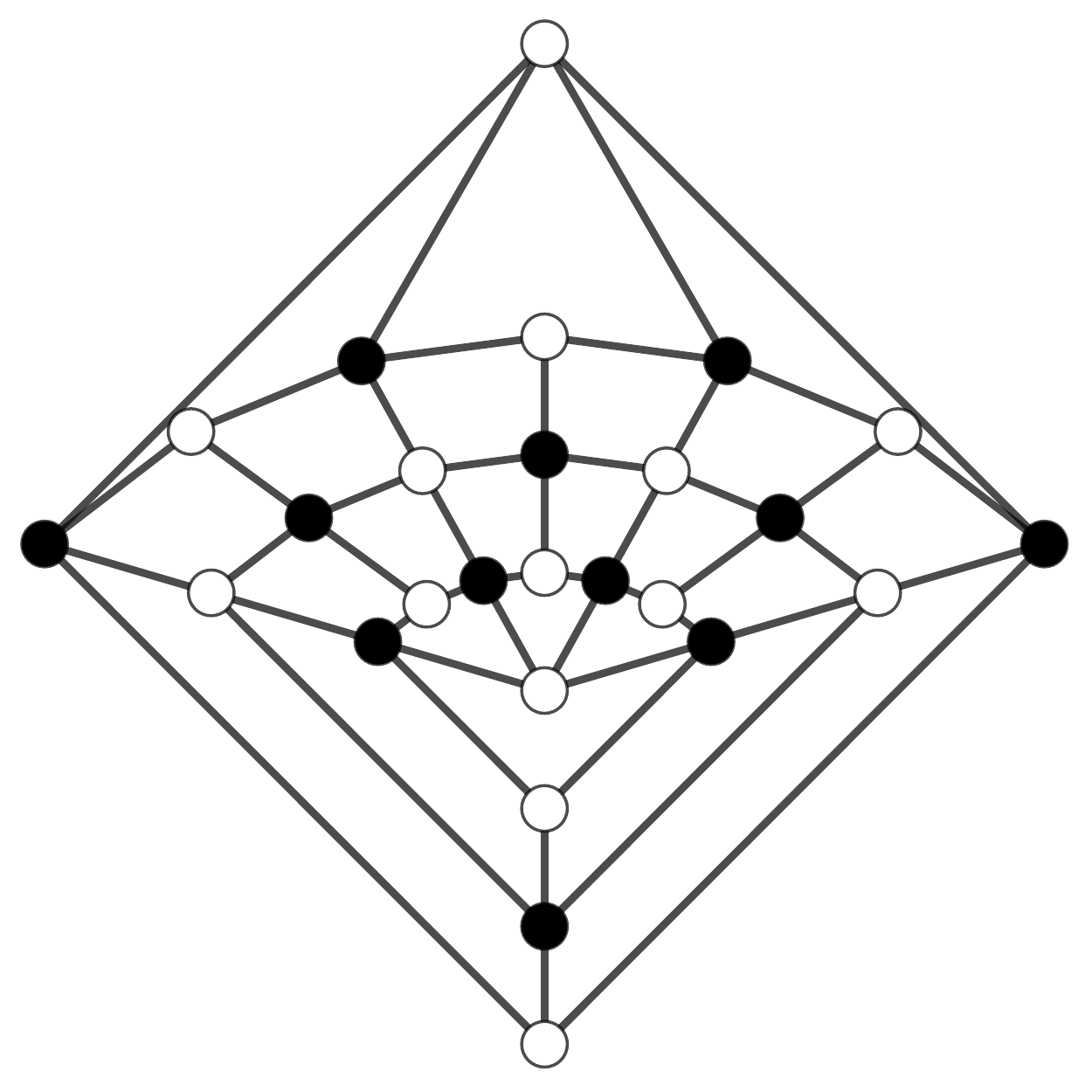}};
\node at (0,-4) {\includegraphics[width=0.31\textwidth]{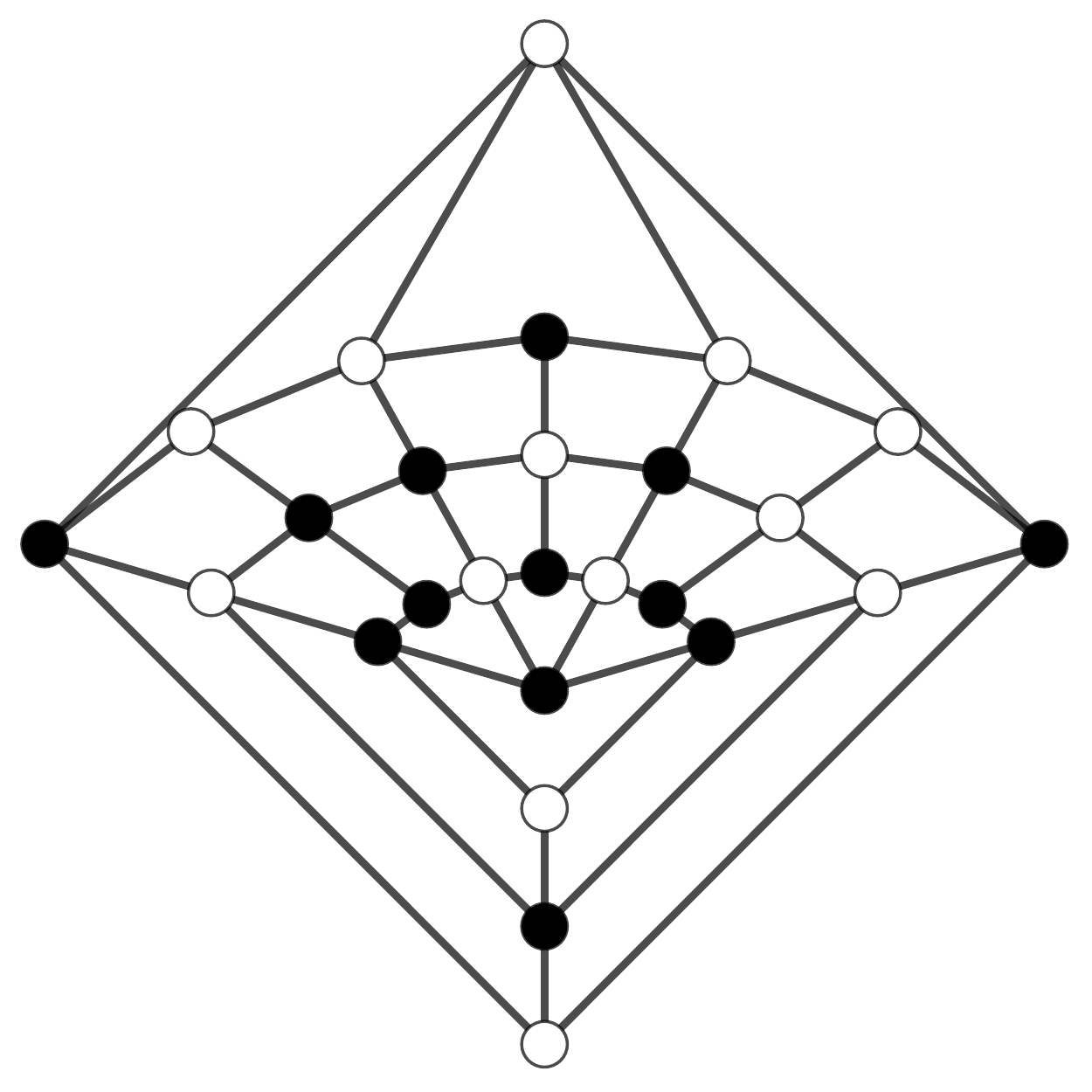}};
\node at (0,-8) {\includegraphics[width=0.31\textwidth]{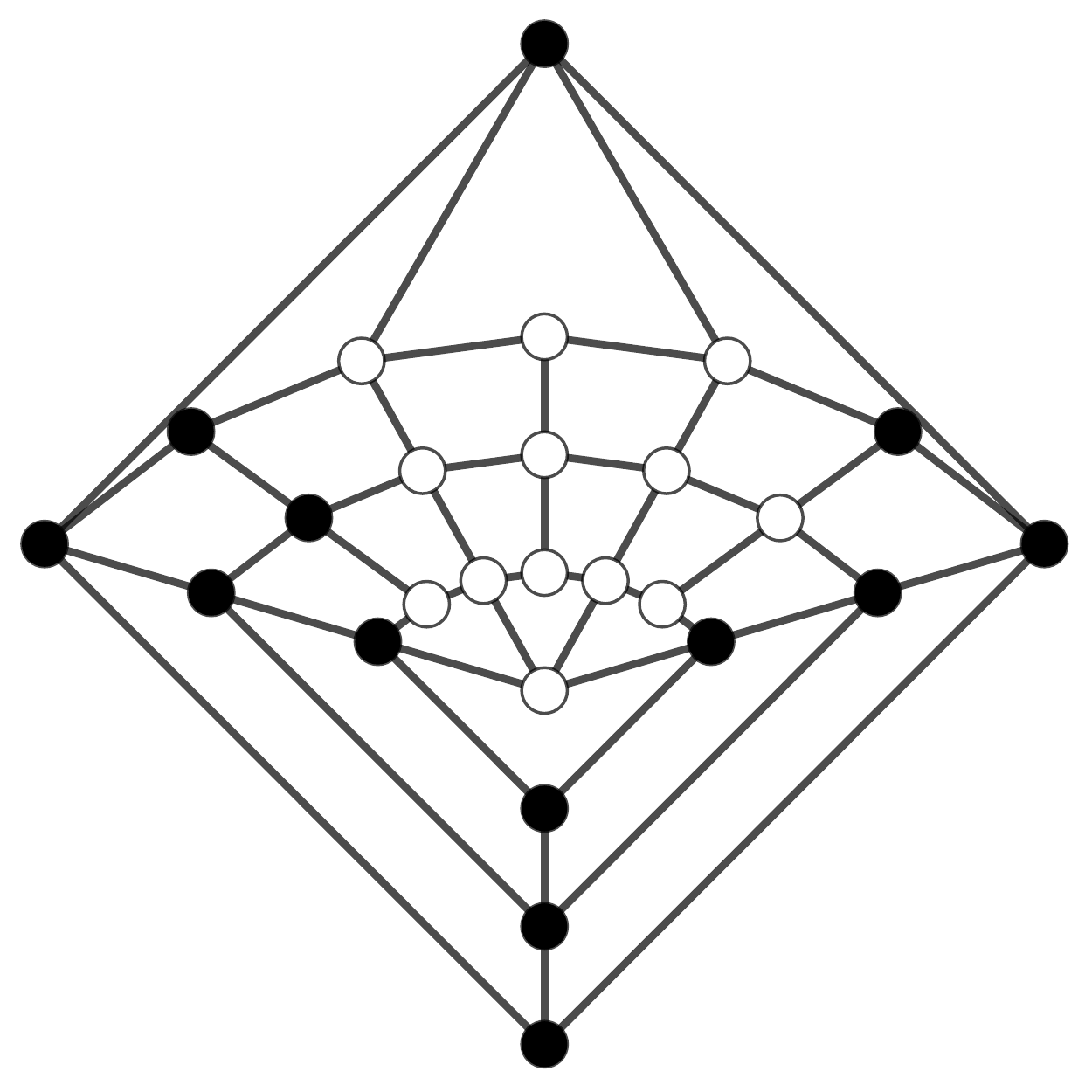}};
\end{tikzpicture}
\vspace{-10pt}
\caption{Sign of the first and second eigenvector on the Deltoidal Icositetrahedron graph as well as sign of their product.}
\label{fig:5}
\end{figure}
\end{center}
\vspace{-5pt}

We quickly discuss this inequality for the cycle graph $C_n$. The eigenvectors are given by, for $ 1\leq k \leq n/2$,
\begin{align*}
x_k(v) = \sin\left( \frac{2 \pi k v}{n} \right) \qquad \mbox{and} \qquad y_k(v) = \cos\left( \frac{2 \pi k v}{n} \right).
\end{align*}
$x_0 \equiv 0$ is ignored. If $n$ is even, then $x_{n/2} \equiv 0$ and is ignored. Eigenvectors $x_k$ and $y_k$ correspond to eigenvalue
$2 - 2 \cos{(2 \pi k/n)}$. Let now $n$ be odd and consider $x_{(n-1)/2}$ and $y_{(n-1)/2}$. The eigenvalue corresponding to these two vectors
is 
$$ 2 - 2\cos\left( \frac{2\pi (n/2-1/2)}{n} \right) = 2 - 2\cos\left( \pi - \frac{\pi}{n} \right) = 4 - \frac{\pi^2}{n^2} +\mbox{l.o.t.}$$
Using the identity $2 \cos{(x)} \sin{(x)} = \sin{(2x)}$, we can write the pointwise product of the eigenvectors as
\begin{align*}
 x_{(n-1)/2}(v) y_{(n-1)/2}(v) &= \frac{1}{2} \sin\left( 2\frac{2 \pi ((n-1)/2) v}{n} \right) = \frac{1}{2} \sin\left( \frac{2 \pi (n-1) v}{n} \right) \\
 &=  \frac{1}{2} \sin\left( \frac{2 \pi (-1) v}{n} \right) = - \frac{y_1(v)}{2} 
\end{align*}
Therefore the quadratic form of the product evaluates to
$$ \left\langle \frac{ x_{(n-1)/2} y_{(n-1)/2}}{\|x_{(n-1)/2} y_{(n-1)/2}\|}, L \frac{ x_{(n-1)/2} y_{(n-1)/2}}{\|x_{(n-1)/2} y_{(n-1)/2}\|} \right\rangle = 2 - 2\cos{\left(\frac{2\pi}{n}\right)} = \frac{4 \pi^2}{n^2} + \mbox{l.o.t}.$$
We can now compare this to the bound obtained in Theorem 1 with $\phi = x_{(n-1)/2}$ and $\psi = y_{(n-1)/2}$. All arising $\ell^{\infty}-$norms are 1. Theorem 1 can be applied with
$ \varepsilon = \pi^2/n^2 + \mbox{l.o.t.}$
which guarantees that
\begin{align*}
 \left\langle \frac{\phi \cdot \psi}{\| \phi \cdot \psi \|_{\ell^2}}, L \frac{ \phi \cdot \psi}{\| \phi \cdot \psi \|_{\ell^2}} \right\rangle &\leq  2\varepsilon \left( \frac{\|\phi\|_{\ell^{\infty}}^2   \|\psi\|_{\ell^2}^2}{\| \phi \cdot \psi \|_{\ell^2}} +  \frac{    \|\phi\|_{\ell^2}^2 \|\psi\|_{\ell^{\infty}}^2 }{\| \phi \cdot \psi \|_{\ell^2}}\right) \\
 &=(1+o(1)) \frac{2 \pi^2}{n^2} \left( \frac{ \|\psi\|_{\ell^2}^2}{\| \phi \cdot \psi \|_{\ell^2}} +  \frac{    \|\phi\|_{\ell^2}^2 }{\| \phi \cdot \psi \|_{\ell^2}}\right) \\
 &= (1+o(1))  \frac{8 \pi^2}{n^2}.
\end{align*}
This shows that the inequality is sharp up to a factor of at most 2.

\subsection{Proof of Theorem 2}
\begin{proof} Let first $i \in V$ be an arbitrary vertex.
We rewrite the equation as
$$ \phi(i) = - \frac{1}{d} \sum_{j \sim_E i} \phi(j) +  \frac{\varepsilon}{d} \phi(i),$$
or, in short, $\phi = D^{-1} \left( \varepsilon \phi - A \phi\right).$ The main idea behind the argument is that it is possible to interpret this as a fixed point equation and iterate it once more. This leads to a cancellation of sign: by then iteratively using the equation in the largest point, it is possible to control the arising error terms. The same idea will also reappear in the proof of the Corollary.
Iterating the equation, we get
\begin{align*}
\phi(i) &= - \frac{1}{d} \sum_{j \sim_E i} \phi(j) +  \frac{\varepsilon}{d} \phi(i) \\
 &=  - \frac{1}{d} \sum_{j \sim_E i}  \left(  - \frac{1}{d} \sum_{k \sim_E j} \phi(k) +  \frac{\varepsilon}{d} \phi(j) \right) +  \frac{\varepsilon}{d} \phi(i) \\
 &=  \frac{1}{d^2} \sum_{j \sim_E i} \sum_{k \sim_E j} \phi(k) -  \frac{\varepsilon}{d^2} \sum_{j \sim_E i}  \phi(j) +  \frac{\varepsilon}{d} \phi(i).
\end{align*}
We note that $(D+A)\phi = \varepsilon \phi$ implies $(D-A)\phi = (2d - \varepsilon) \phi$ and thus
$$ -  \frac{\varepsilon}{d^2} \sum_{j \sim_E i}  \phi(j) +  \frac{\varepsilon}{d} \phi(i) = \frac{\varepsilon}{d} \left( - AD^{-1} \phi + \phi\right)(i) = \frac{\varepsilon}{d} \left(2 - \frac{\varepsilon}{d}\right) \phi(i).$$
Let us now assume that $\|\phi\|_{\ell^{\infty}(V)}$ is assumed in a vertex $m \in V$, i.e. $\phi(m) = \| \phi\|_{\ell^{\infty}(V)}$ (this can always be
assumed after potentially replacing $\phi$ by $-\phi$). Then, for all $i \in V$, we have the identity
$$
\left(1 - \frac{\varepsilon}{d} \right)^2 \phi(i) =   \frac{1}{d^2} \sum_{j \sim_E i}  \sum_{k \sim_E j} \phi(k).
$$
Therefore
$$  \sum_{j \in V} \left[(AD^{-1})^{2k} \delta_m\right](j) \cdot \phi(j) = \left(1 - \frac{ \varepsilon}{d} \right)^{2k} \phi(m) .$$
Using Cauchy-Schwarz, we deduce the desired result since
$$  \left(1 - \frac{ \varepsilon}{d} \right)^{2k} |\phi(m)| \leq  \left(\max_{m \in V} \sum_{j =1}^{n} ((AD^{-1})^{2k})_{m, j}^2 \right)^{1/2} \| \phi\|_{\ell^2}.$$
\end{proof}

\vspace{-10pt}
\begin{center}
\begin{figure}[h!]
\begin{tikzpicture}[rotate=90]
\node at (0,0) {\includegraphics[width=0.29\textwidth]{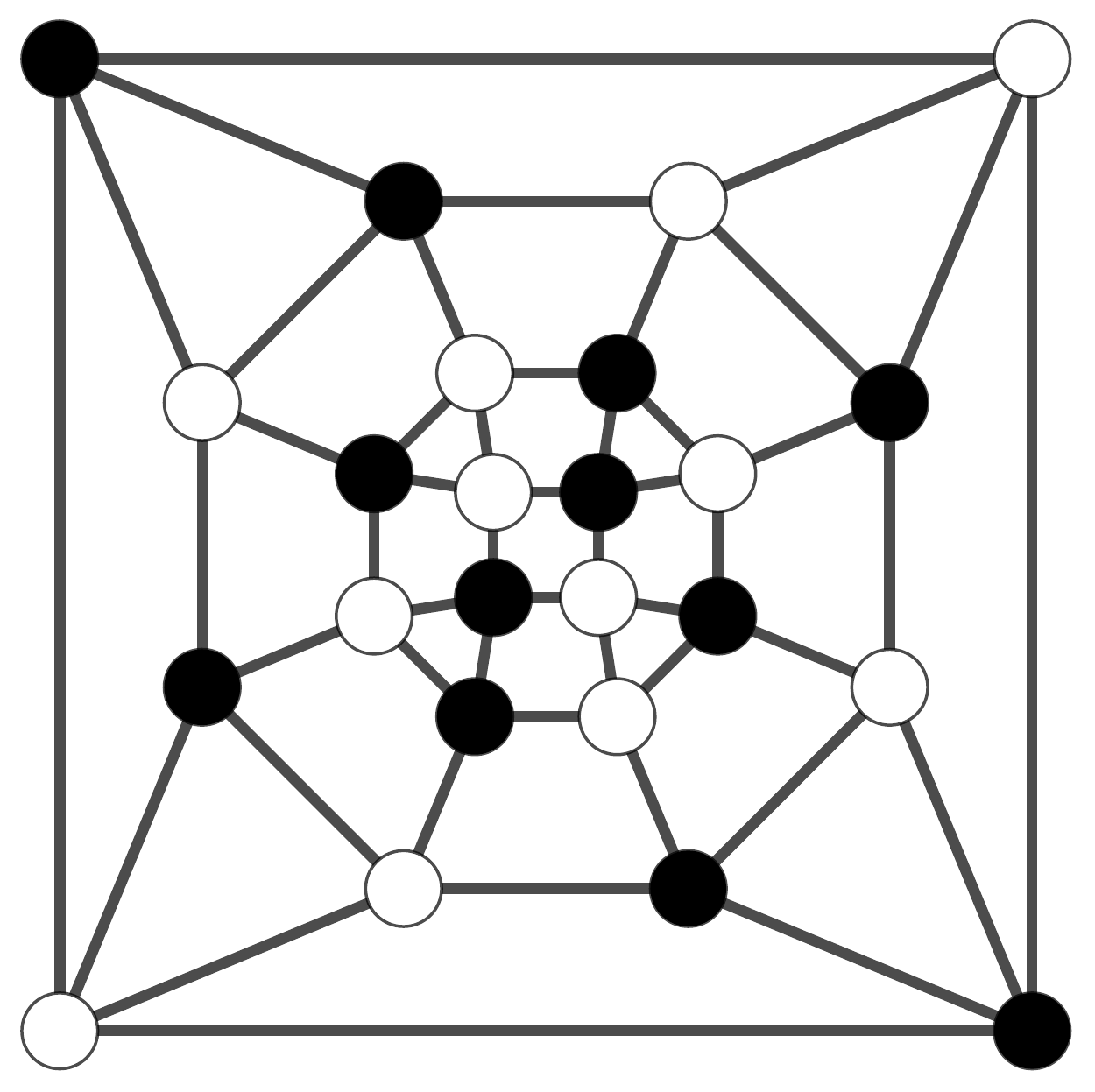}};
\node at (0,-4) {\includegraphics[width=0.29\textwidth]{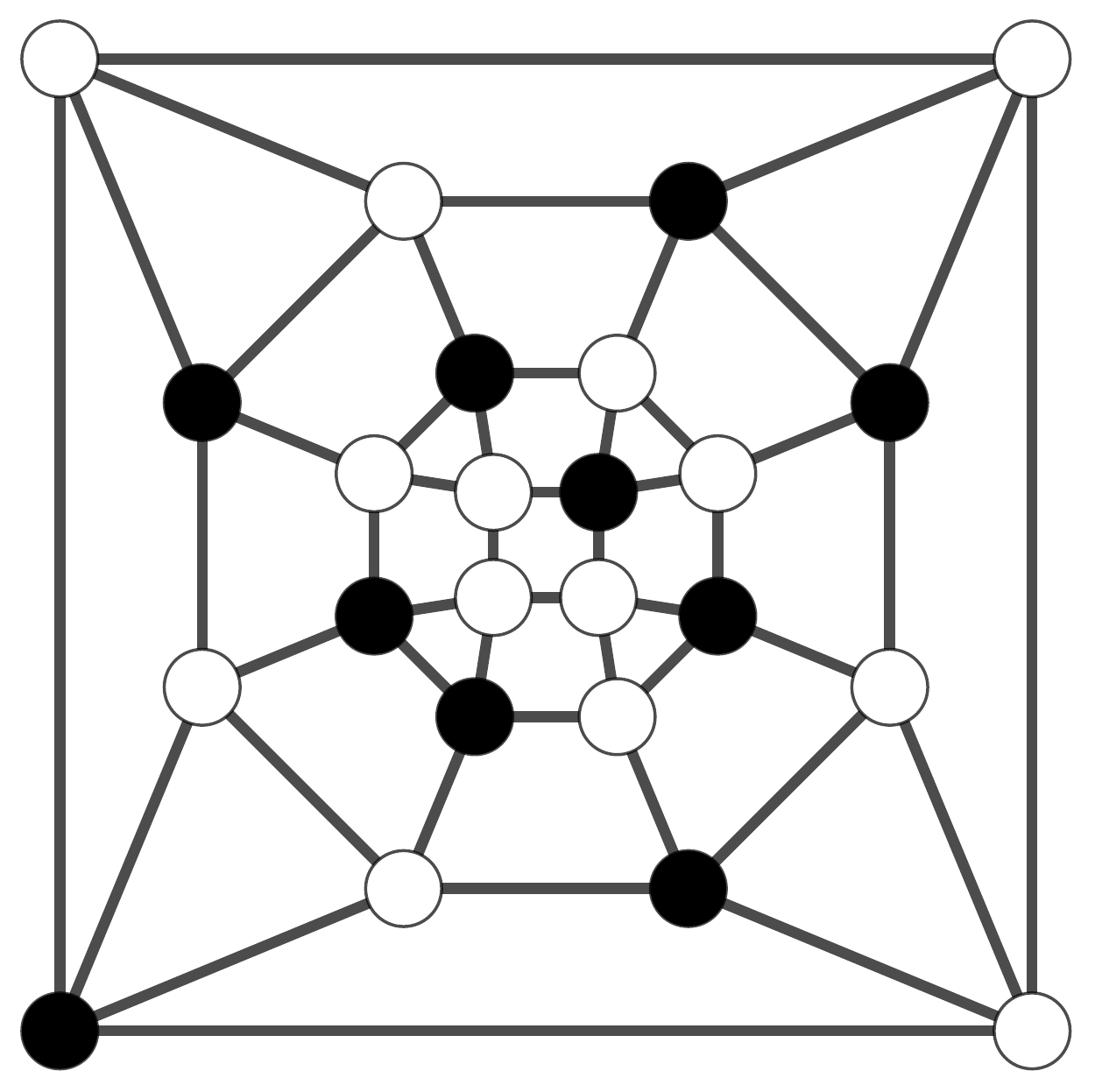}};
\node at (0,-8) {\includegraphics[width=0.29\textwidth]{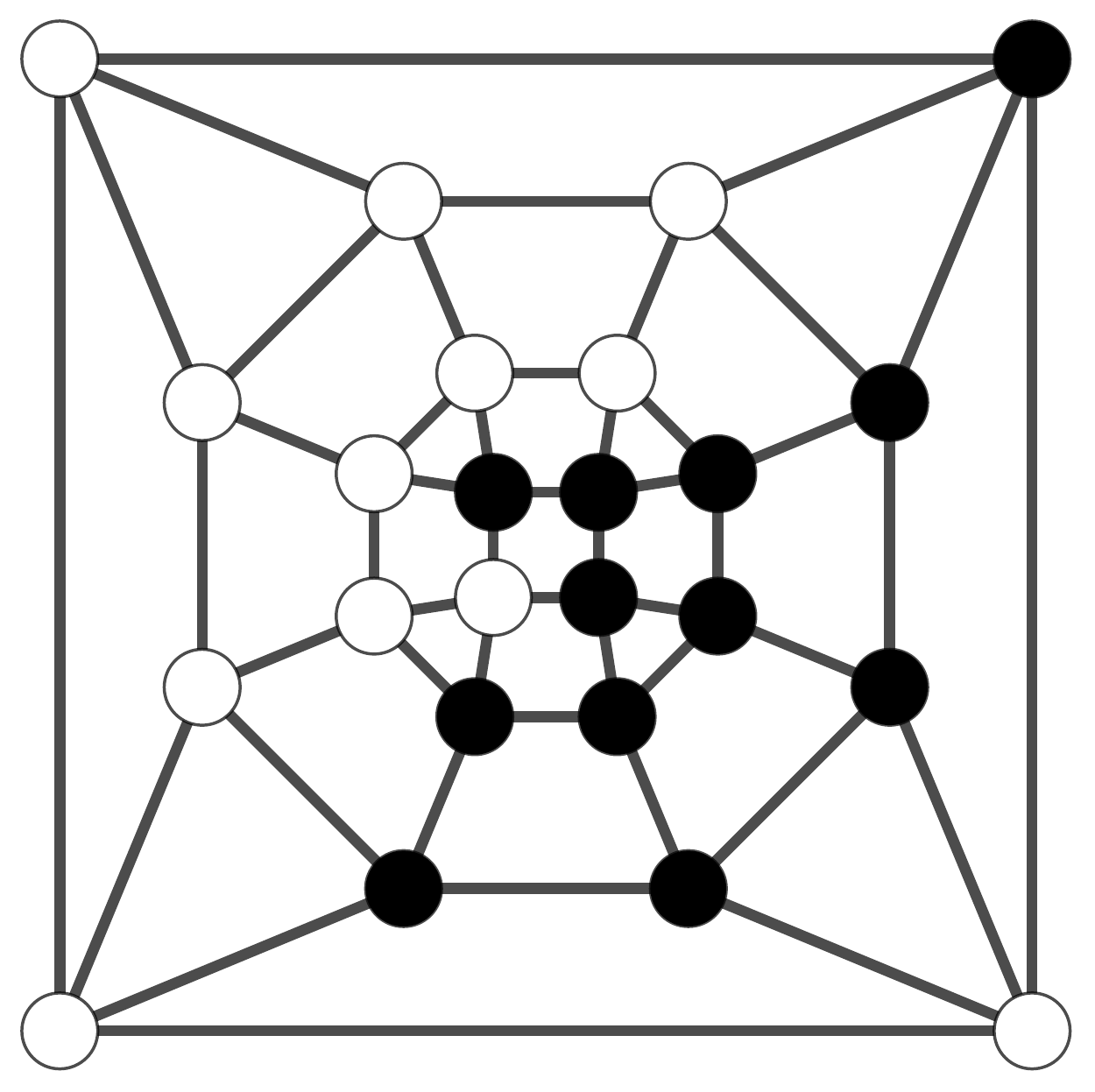}};
\end{tikzpicture}
\vspace{-10pt}
\caption{Sign of the first and third eigenvector on the small rhombicuboctahedral graph as well as sign of their product.}
\label{fig:9}
\end{figure}
\end{center}
\vspace{-10pt}

\subsection{Proof of the Corollary}
\begin{proof} The Corollary follows quickly from considering the same approach as in the proof of Theorem 2.
The equation
$$ (d \cdot \mbox{Id}_{n \times n}+A)\phi = \varepsilon \phi$$
can be rewritten as
$$ \phi = \frac{1}{d} \left( \varepsilon \phi - A \phi \right).$$
Plugging this identity into itself, we obtain
\begin{align*}
 \phi &= \frac{1}{d} \left( \varepsilon \phi - A \left[\frac{1}{d} \left( \varepsilon \phi - A \phi \right) \right] \right) \\
 &= \frac{1}{d^2} A^2 \phi -  \frac{\varepsilon}{d^2} A \phi + \frac{\varepsilon}{d} \phi \\
 &=    \frac{\varepsilon}{d^2} A \phi -  \frac{\varepsilon}{d^2} A \phi + \frac{\varepsilon}{d} \phi  
 \end{align*}
 Noting, as above, that
 $$ -  \frac{\varepsilon}{d^2} A \phi + \frac{\varepsilon}{d} \phi   = - \frac{\varepsilon}{d} \left( AD^{-1} \phi - \phi\right) = \frac{\varepsilon}{d}\left(2 - \frac{\varepsilon}{d}\right) \phi$$
 then yields the desired result.
\end{proof}

\end{document}